# Divergence-free $H(\mathrm{div})$-FEM for time-dependent incompressible flows with applications to high Reynolds number vortex dynamics

**Philipp W. Schroeder** · **Gert Lube**

**Abstract** In this article, we consider exactly divergence-free $H(\mathrm{div})$-conforming finite element methods for time-dependent incompressible viscous flow problems. This is an extension of previous research concerning divergence-free $H^1$-conforming methods. For the linearised Oseen case, the first semi-discrete numerical analysis for time-dependent flows is presented whereby special emphasis is put on pressure- and Reynolds-semi-robustness. For convection-dominated problems, the proposed method relies on a velocity jump upwind stabilisation which is not gradient-based. Complementing the theoretical results, $H(\mathrm{div})$-FEM are applied to the simulation of full nonlinear Navier–Stokes problems. Focussing on dynamic high Reynolds number examples with vortical structures, the proposed method proves to be capable of reliably handling the planar lattice flow problem, Kelvin–Helmholtz instabilities and freely decaying two-dimensional turbulence.

## 1 Introduction

In this paper, we consider time-dependent incompressible flows fulfilling [46,41,20]

$$\begin{cases} \partial_t \boldsymbol{u} - \nu\Delta\boldsymbol{u} + (\boldsymbol{\beta}\cdot\nabla)\boldsymbol{u} + \nabla p = \boldsymbol{f} & \text{in } (0,T]\times\Omega, \quad (1a)\\ \nabla\cdot\boldsymbol{u} = 0 & \text{in } (0,T]\times\Omega, \quad (1b)\\ \boldsymbol{u} = \boldsymbol{0} & \text{on } [0,T]\times\partial\Omega, \quad (1c)\\ \boldsymbol{u}(0,\boldsymbol{x}) = \boldsymbol{u}_0(\boldsymbol{x}) & \text{for } \boldsymbol{x}\in\Omega. \quad (1d) \end{cases}$$

Philipp W. Schroeder & Gert Lube
Institute for Numerical and Applied Mathematics
Georg-August-University Göttingen
D-37083 Göttingen, Germany
E-mail: {p.schroeder,lube}@math.uni-goettingen.de
ORCID P.W. Schroeder: https://orcid.org/0000-0001-7644-4693

For the space dimension $d \in \{2,3\}$, $\Omega \subset \mathbb{R}^d$ denotes a bounded connected domain. Moreover, $\boldsymbol{u}\colon (0,T] \times \Omega \to \mathbb{R}^d$ represents the velocity field. For $\boldsymbol{\beta} = \boldsymbol{u}$, the nonlinear Navier–Stokes equations are recovered. A linearisation of the problem, known as the Oseen problem, can be considered whenever an *a priori* known convective velocity $\boldsymbol{\beta}\colon (0,T] \times \Omega \to \mathbb{R}^d$ is taken instead. In both cases, $p\colon (0,T] \times \Omega \to \mathbb{R}$ is the (zero-mean) kinematic pressure, $\boldsymbol{f}\colon (0,T] \times \Omega \to \mathbb{R}^d$ represents external body forces and $\boldsymbol{u}_0\colon \Omega \to \mathbb{R}^d$ stands for a suitable initial condition for the velocity. The underlying fluid is assumed to be Newtonian with constant (dimensionless) kinematic viscosity $0 < \nu \ll 1$.

We want to use an $\boldsymbol{H}(\mathrm{div})$-conforming, inf-sup stable and exactly divergence-free finite element method (FEM) for solving (1) approximately. Divergence-free FEMs enable a strict separation between the approximation of velocity and pressure in the sense that it is possible to obtain error estimates where the quality of the pressure approximation does not influence the velocity error. This property is referred to as 'pressure-robustnes'; cf. [32,35].

Moreover, using an $\boldsymbol{H}(\mathrm{div})$-conforming FEM allows the usage of the whole machinery known from Discontinuous Galerkin FEM (dG-FEM). Especially an upwind treatment of the convective term can be incorporated quite naturally. In this context, whenever problems with high Reynolds numbers *Re* (equivalently with small viscosity $\nu$) are considered, it is sensible to strive for methods for which '*Re*-semi-robustness' can be shown—that is, methods whose error estimates, including Gronwall constants, do not explicitly depend on *Re* (equivalently on $\nu^{-1}$); cf. [39]. The price to be paid usually are certain regularity assumptions for the exact solution.

Thus, this work can be seen as an advancement and extension of the authors' previous work [42], where divergence-free $\boldsymbol{H}^1$-conforming FEM have been analysed for the time-dependent Navier–Stokes problem. For a literature overview concerning $\boldsymbol{H}^1$-conforming FEM, we also refer to [42]. The present contribution is split into the following parts:

- **Numerical analysis:** For a semi-discrete numerical error analysis, we consider the time-dependent Oseen equations with a known convective field $\boldsymbol{\beta}$. The analysis of linearised incompressible flow problems shall act as a proof of concept and a first step towards the nonlinear Navier–Stokes problem. Moreover, due to the above mentioned possibility of strictly separating velocity and pressure, this work exclusively focusses on the velocity approximation.
- **Application:** In contrast to the theoretical part, for our numerical experiments, we apply the $\boldsymbol{H}(\mathrm{div})$-FEM to the full nonlinear Navier–Stokes problem with $\boldsymbol{\beta} = \boldsymbol{u}$. There, for high Reynolds number problems, we assess the quality of $\boldsymbol{H}(\mathrm{div})$-FEM in comparison with other discretisation schemes. The considered problems all have a strong vortical structure and vary in their behaviour from rather static to considerably dynamic.

Concerning previous research, in [49,48], $\boldsymbol{H}(\mathrm{div})$-conforming and divergence-free FEMs are introduced for the stationary Stokes problem; but they do not provide pressure-robust error estimates. On the other hand, in [32], a pressure-robust estimate for $\boldsymbol{H}(\mathrm{div})$-FEM in a discrete energy norm is given for the stationary Stokes problem; but without many details.

For the stationary Navier–Stokes equations, $\boldsymbol{H}(\mathrm{div})$-FEM are considered in [16]. However, for the error analysis, one is referred to the Local Discontinuous Galerkin method [15] and in this work, neither pressure- nor *Re*-semi-robustness plays a role. Actually, in the context of stationary problems, the high Reynolds number case is usually excluded. For the time-dependent incompressible Euler equations ($\nu = 0$), on the other hand, $\boldsymbol{H}(\mathrm{div})$-FEM are applied in [27]. But viscous effects do not play a role in the Euler problem. Therefore, to the best of the authors' knowledge, this contribution is the first one which analyses and applies $\boldsymbol{H}(\mathrm{div})$-FEM to *time-dependent viscous* incompressible flows. Additionally, we attach great importance to the aspects of pressure- and *Re*-semi-robustness.

*Organisation of the article:* In Section 2, we introduce divergence-free $\boldsymbol{H}(\mathrm{div})$-conforming and inf-sup stable FEM for time-dependent incompressible flow problems. This includes weak formulations with their corresponding function spaces, assumptions which are needed subsequently and a brief discussion of the stability and well-posedness of the proposed method. Then, Section 3 is concerned with the derivation of *a priori*, pressure- and *Re*-semi-robust error estimates for the spatially discretised velocity of the linearised problem. Finally, in Section 4, the $\boldsymbol{H}(\mathrm{div})$-FEM is applied to the simulation of high Reynolds number flows for the full nonlinear Navier–Stokes problem. There, we present numerical examples for the planar lattice flow problem, Kelvin–Helmholtz instabilities and freely decaying two-dimensional turbulence.

## 2 $\boldsymbol{H}$(div)-FEM for time-dependent incompressible flows

*Notation:* In what follows, for $K \subseteq \Omega$ we use the standard Sobolev spaces $W^{m,p}(K)$ for scalar-valued functions with associated norms $\|\cdot\|_{W^{m,p}(K)}$ and seminorms $|\cdot|_{W^{m,p}(K)}$ for $m \geqslant 0$ and $p \geqslant 1$. Spaces and norms for vector- and tensor-valued functions are indicated with bold letters. We obtain the Lebesgue space $W^{0,p}(K) = L^p(K)$ and the Hilbert space $W^{m,2}(K) = H^m(K)$. Additionally, the closed subspaces $H_0^1(K)$ consisting of $H^1(K)$-functions with vanishing trace on $\partial K$ and the set $L_0^2(K)$ of $L^2(K)$-functions with zero mean in $K$ play an important role. The $L^2(K)$-inner product is denoted by $(\cdot,\cdot)_K$ and, if $K = \Omega$, we usually omit the domain completely when no confusion can arise. Furthermore, with regard to time-dependent problems, given a Banach space $\boldsymbol{X}$ and a time instance $t^*$, the Bochner space $L^p(0,t^*;\boldsymbol{X})$ for $p \in [1,\infty]$ is used. In the case $t^* = T$, we frequently use the abbreviation $L^p(\boldsymbol{X}) = L^p(0,T;\boldsymbol{X})$.

### 2.1 Time-dependent Oseen problem

With $\boldsymbol{V} = \boldsymbol{H}_0^1(\Omega)$ and $Q = L_0^2(\Omega)$, we introduce the spaces for velocity and pressure as

$$\boldsymbol{V}^T = \left\{ \boldsymbol{v} \in L^2(0,T;\boldsymbol{V}) \colon \partial_t \boldsymbol{v} \in L^2\big(0,T;\boldsymbol{L}^2\big) \right\}, \qquad Q^T = L^2(0,T;Q). \tag{2}$$

Note that throughout this work, we thus assume at least the mild regularity $\partial_t \boldsymbol{u} \in L^2\big(0,T;\boldsymbol{L}^2\big)$ for the exact velocity. Then, provided the forcing term $\boldsymbol{f}$ is sufficiently smooth, the following well-known variational formulation of problem (1) on the continuous level is obtained:

$$\begin{cases} ?\,(\boldsymbol{u},p) \in \boldsymbol{V}^T \times Q^T \text{ with } \boldsymbol{u}(0) = \boldsymbol{u}_0 \text{ s.t., } \forall\,(\boldsymbol{v},q) \in \boldsymbol{V} \times Q, & \text{(3a)} \\ (\partial_t \boldsymbol{u},\boldsymbol{v}) + \nu a(\boldsymbol{u},\boldsymbol{v}) + c(\boldsymbol{\beta};\boldsymbol{u},\boldsymbol{v}) + b(\boldsymbol{v},p) - b(\boldsymbol{u},q) = (\boldsymbol{f},\boldsymbol{v}). & \text{(3b)} \end{cases}$$

The multilinear forms are given by

$$a(\boldsymbol{w},\boldsymbol{v}) = \int_\Omega \nabla \boldsymbol{w} : \nabla \boldsymbol{v} \,\mathrm{d}\boldsymbol{x}, \qquad c(\boldsymbol{\beta};\boldsymbol{w},\boldsymbol{v}) = \int_\Omega (\boldsymbol{\beta}\cdot\nabla)\boldsymbol{w}\cdot\boldsymbol{v}\,\mathrm{d}\boldsymbol{x}, \tag{4a}$$

$$b(\boldsymbol{w},q) = -\int_\Omega q(\nabla\cdot\boldsymbol{w})\,\mathrm{d}\boldsymbol{x}. \tag{4b}$$

In this work, $\boldsymbol{\beta}\colon (0,T]\times\Omega\to\mathbb{R}^d$ denotes the known convective velocity. We assume that $\boldsymbol{\beta}\in L^\infty(0,T;\boldsymbol{L}^\infty)$ with $\nabla\cdot\boldsymbol{\beta}=0$ pointwise and $\boldsymbol{\beta}\cdot\boldsymbol{n}\big|_{\partial\Omega}=0$, where $\boldsymbol{n}$ denotes the outer unit normal vector to $\partial\Omega$. In applications, the field $\boldsymbol{\beta}$ can be thought of as an approximation of $\boldsymbol{u}$. We abbreviate $\|\boldsymbol{\beta}\|_{L^\infty(\boldsymbol{L}^\infty)}=\|\boldsymbol{\beta}\|_\infty$. Weakly divergence-free velocities belong to

$$\boldsymbol{V}^{\mathrm{div}}=\{\boldsymbol{v}\in\boldsymbol{V}\colon b(\boldsymbol{v},q)=0,\ \forall q\in Q\}. \tag{5}$$

### 2.2 Discrete setting and assumptions

In this work, we focus on FEM which are $\boldsymbol{H}(\mathrm{div})$-conforming, where

$$\boldsymbol{H}(\mathrm{div};\Omega)=\left\{\boldsymbol{w}\in\boldsymbol{L}^2(\Omega)\colon\ \nabla\cdot\boldsymbol{w}\in L^2(\Omega)\right\}. \tag{6}$$

Let $\mathcal{T}_h$ be a shape-regular FE partition of $\Omega$ without hanging nodes and mesh size $h=\max_{K\in\mathcal{T}_h}h_K$, where $h_K$ denotes the diameter of the particular element $K\in\mathcal{T}_h$. Since the subsequent velocity approximation will not be $\boldsymbol{H}^1$-conforming, the broken Sobolev space

$$\boldsymbol{H}^m(\mathcal{T}_h)=\left\{\boldsymbol{w}\in\boldsymbol{L}^2(\Omega)\colon \boldsymbol{w}\big|_K\in\boldsymbol{H}^m(K),\ \forall K\in\mathcal{T}_h\right\} \tag{7}$$

is introduced. Define the broken gradient $\nabla_h\colon \boldsymbol{H}^1(\mathcal{T}_h)\to\boldsymbol{L}^2(\Omega)$ by $(\nabla_h\boldsymbol{w})\big|_K=\nabla\big(\boldsymbol{w}\big|_K\big)$.

The skeleton $\mathcal{F}_h$ denotes the set of all facets with $\mathcal{F}_K=\{F\in\mathcal{F}_h\colon F\subset\partial K\}$ and $N_\partial=\max_{K\in\mathcal{T}_h}\mathrm{card}(\mathcal{F}_K)$. Moreover, $\mathcal{F}_h=\mathcal{F}_h^i\cup\mathcal{F}_h^\partial$ where $\mathcal{F}_h^i$ is the subset of interior facets and $\mathcal{F}_h^\partial$ collects all boundary facets $F\subset\partial\Omega$. To any $F\in\mathcal{F}_h$ we assign a unit normal vector $\boldsymbol{n}_F$ where, for $F\in\mathcal{F}_h^\partial$, this is the outer unit normal vector $\boldsymbol{n}$. If $F\in\mathcal{F}_h^i$, there are two adjacent elements $K^+$ and $K^-$ sharing the facet $F=\overline{\partial K^+}\cap\overline{\partial K^-}$ and $\boldsymbol{n}_F$ points in an arbitrary but fixed direction. Let $\phi$ be any piecewise smooth (scalar-, vector- or tensor-valued) function with traces from within the interior of $K^\pm$ denoted by $\phi^\pm$, respectively. Then, we define the jump $[\![\cdot]\!]_F$ and average $\{\!\!\{\cdot\}\!\!\}_F$ operator across interior facets $F\in\mathcal{F}_h^i$ by

$$[\![\phi]\!]_F=\phi^+-\phi^-\quad\text{and}\quad\{\!\!\{\phi\}\!\!\}_F=\frac{1}{2}\big(\phi^++\phi^-\big). \tag{8}$$

For boundary facets $F\in\mathcal{F}_h^\partial$ we set $[\![\phi]\!]_F=\{\!\!\{\phi\}\!\!\}_F=\phi$. These operators act componentwise for vector- and tensor-valued functions. Frequently, the subscript indicating the facet is omitted. It is important to have in mind the following characterisation of $\boldsymbol{H}(\mathrm{div})$-functions.

**LEMMA 2.1 (*Characterisation of $\boldsymbol{H}(\mathrm{div};\Omega)$*)**

*Let $\boldsymbol{w}\in\boldsymbol{H}^1(\mathcal{T}_h)$. If $[\![\boldsymbol{w}]\!]\cdot\boldsymbol{n}_F=0$ for all $F\in\mathcal{F}_h^i$, then $\boldsymbol{w}\in\boldsymbol{H}(\mathrm{div};\Omega)$.*

**PROOF:** Cf., for example, [19, Lemma 1.24]. ■

In the following, $\mathbb{P}_k(K)$ (vector-valued: $\boldsymbol{\mathbb{P}}_k(K)$) denotes either the space of all polynomials on $K$ with degree less or equal to $k$ (simplicial mesh) or with degree less or equal to $k$ in each variable (tensor-product mesh), depending on the particular situation.

Let us turn to specifying the assumptions which are needed in the remainder of this paper. Define the following discrete FE spaces for velocity and pressure, respectively:

$$\boldsymbol{V}_h=\left\{\boldsymbol{v}_h\in\boldsymbol{H}(\mathrm{div};\Omega)\colon \boldsymbol{v}_h\big|_K\in\boldsymbol{V}_k(K),\ \forall K\in\mathcal{T}_h;\ \boldsymbol{v}_h\cdot\boldsymbol{n}\big|_{\partial\Omega}=0\right\}\not\subset\boldsymbol{V} \tag{9a}$$

$$Q_h=\left\{q_h\in L^2_0(\Omega)\colon q_h\big|_K\in\mathbb{P}_\ell(K),\ \forall K\in\mathcal{T}_h\right\}\subset Q \tag{9b}$$

Here, the local space $\boldsymbol{V}_k(K)$ is a set of vector-valued piecewise polynomials of order $k \geqslant 1$, which, in order to keep the theory in Section 3 widely applicable, is not explicitly specified further. Instead, we shall only introduce several required global and local properties for the resulting FE pair. For the pressure space, $\ell \in \{k-1, k\}$ is assumed. Some valid and explicit examples for $\boldsymbol{V}_h/Q_h$ are given below.

**Assumption A1:** The global spaces $\boldsymbol{V}_h$ and $Q_h$ are divergence-conforming. That is,

$$\nabla\cdot\boldsymbol{V}_h \subseteq Q_h. \tag{10}$$

Property (10) ensures that the velocity approximation will be exactly divergence-free [32].

**Assumption A2:** The global spaces $\boldsymbol{V}_h$ and $Q_h$ form a discretely inf-sup stable FE pair. That is, there exists $\beta > 0$, independent of the mesh size $h$, such that

$$\inf_{q_h\in Q_h\setminus\{0\}} \sup_{\boldsymbol{v}_h\in\boldsymbol{V}_h\setminus\{\boldsymbol{0}\}} \frac{b(\boldsymbol{v}_h,q_h)}{|||\boldsymbol{v}_h|||_e \,\|q_h\|_{L^2}} \geqslant \beta. \tag{11}$$

Here, $|||\cdot|||_e$ denotes a suitable energy norm. Due to the $\boldsymbol{H}(\mathrm{div})$-conformity of $\boldsymbol{V}_h$, the pressure-velocity coupling $b(\cdot,\cdot)$ does not have to be modified in the discrete setting. Note that (11) ensures that the space of discretely divergence-free velocities $\boldsymbol{V}_h^{\mathrm{div}}$ is non-trivial, that is

$$\boldsymbol{V}_h^{\mathrm{div}} = \{\boldsymbol{v}_h \in \boldsymbol{V}_h \colon b(\boldsymbol{v}_h,q_h) = 0,\ \forall q_h \in Q_h\} \neq \{\boldsymbol{0}\}. \tag{12}$$

**Assumption A3:** The space $\boldsymbol{V}_h$ has optimal approximation properties in the following sense. There exists a velocity approximation operator $\boldsymbol{j}_h \colon \boldsymbol{V} \to \boldsymbol{V}_h$ such that, for all $\boldsymbol{w} \in \boldsymbol{H}^r(\Omega)$ with $r > 3/2$ and $r_{\boldsymbol{u}} = \min\{r, k+1\}$,

$$\|\boldsymbol{w}-\boldsymbol{j}_h\boldsymbol{w}\|_{\boldsymbol{L}^2(K)} + h_K\,\|\boldsymbol{w}-\boldsymbol{j}_h\boldsymbol{w}\|_{\boldsymbol{H}^1(K)} \leqslant C h_K^{r_{\boldsymbol{u}}} |\boldsymbol{w}|_{\boldsymbol{H}^{r_{\boldsymbol{u}}}(K)}, \quad \forall K \in \mathcal{T}_h. \tag{13}$$

A direct consequence of the optimal approximation property (13), together with a continuous trace inequality [19], is the ability to bound polynomial approximation errors on facets:

$$\|\boldsymbol{w}-\boldsymbol{j}_h\boldsymbol{w}\|_{\boldsymbol{L}^2(F)} + h_K\,\|\nabla(\boldsymbol{w}-\boldsymbol{j}_h\boldsymbol{w})\cdot\boldsymbol{n}_K\|_{\boldsymbol{L}^2(F)} \leqslant C h_K^{r_{\boldsymbol{u}}-\frac{1}{2}} |\boldsymbol{w}|_{\boldsymbol{H}^{r_{\boldsymbol{u}}}(K)}, \quad \forall F \in \mathcal{F}_K,\ K \in \mathcal{T}_h. \tag{14}$$

Furthermore, concerning the pressure, it is well-known that for all $q \in Q \cap H^s(\Omega)$ with $s \geqslant 1$ and $r_p = \min\{s, \ell+1\}$ the local orthogonal $L^2$-projection $\pi_0 \colon L^2(K) \to \mathbb{P}_\ell(K)$ fulfils

$$\|q-\pi_0 q\|_{L^2(K)} \leqslant C h_K^{r_p} |q|_{H^{r_p}(K)}, \quad \forall K \in \mathcal{T}_h. \tag{15}$$

**Assumption A4:** $\boldsymbol{j}_h$ fulfils the following commuting diagram property:

$$\nabla\cdot(\boldsymbol{j}_h\boldsymbol{w}) = \pi_0(\nabla\cdot\boldsymbol{w}) \tag{16}$$

Note that $\nabla\cdot\boldsymbol{w} = 0$ pointwise implies that $\nabla\cdot(\boldsymbol{j}_h\boldsymbol{w}) = 0$ also holds in a pointwise sense. Therefore, property (16) will ensure that our analysis yields pressure-robust error estimates for the velocity. However, note that not every choice of velocity approximation operator automatically leads to pressure-robust estimates.

**Assumption A5:** Let $0 \leqslant m \leqslant \ell$ and $1 \leqslant p, q \leqslant \infty$. The local space $\boldsymbol{V}_k(K)$ satisfies the local inverse inequality [21, Lemma 1.138]

$$\forall \boldsymbol{v}_h \in \boldsymbol{V}_k(K)\colon \quad \|\boldsymbol{v}_h\|_{\boldsymbol{W}^{\ell,p}(K)} \leqslant C_{\mathrm{inv}} h_K^{m-\ell+d\left(\frac{1}{p}-\frac{1}{q}\right)} \|\boldsymbol{v}_h\|_{\boldsymbol{W}^{m,q}(K)}, \quad \forall K \in \mathcal{T}_h. \tag{17}$$

We will not distinguish between different applications of (17) and in the end, $C_{\mathrm{inv}}$ can be thought of as a maximum over all occurring constants. Also on shape-regular meshes, the following assumption is fulfilled.

**Assumption A6:** The space $\boldsymbol{V}_k(K)$ satisfies the discrete trace inequality [19, Remark 1.47]

$$\forall \boldsymbol{v}_h \in \boldsymbol{V}_k(K)\colon \quad \|\boldsymbol{v}_h\|_{\boldsymbol{L}^2(\partial K)} \leqslant C_{\mathrm{tr}} N_\partial^{1/2} h_K^{-1/2} \|\boldsymbol{v}_h\|_{\boldsymbol{L}^2(K)}, \quad \forall K \in \mathcal{T}_h. \tag{18}$$

**Remark 2.2:** Several classical examples of spaces fulfilling these properties can be found in the monograph [7]. For $\ell = k$, let us specifically mention the family of Raviart–Thomas elements on simplicial meshes. We will use this element later for our numerical experiments. For $\ell = k-1$, the family of Brezzi–Douglas–Marini elements on either simplicial or tensor-product meshes is also applicable. Note that both type of elements come with an interpolation operator fulfilling Assumptions 3 and 4. More recent developments in the direction of creating new $\boldsymbol{H}(\mathrm{div})$ elements can be found, for example, in [49]. ▲

### 2.3 FEM and well-posedness

Our discrete space-time velocity and pressure spaces are

$$\boldsymbol{V}_h^T = \left\{\boldsymbol{v}_h \in L^2(0,T;\boldsymbol{V}_h)\colon \partial_t \boldsymbol{v}_h \in L^2(0,T;\boldsymbol{V}_h)\right\}, \quad Q_h^T = L^2(0,T;Q_h). \tag{19}$$

The space-semidiscrete variational formulation of (3) reads as follows:

$$\begin{cases} ?\,(\boldsymbol{u}_h, p_h) \in \boldsymbol{V}_h^T \times Q_h^T \text{ with } \boldsymbol{u}_h(0) = \boldsymbol{u}_{0h} \text{ s.t., } \forall (\boldsymbol{v}_h, q_h) \in \boldsymbol{V}_h \times Q_h, & \text{(20a)} \\ (\partial_t \boldsymbol{u}_h, \boldsymbol{v}_h) + \nu a_h(\boldsymbol{u}_h, \boldsymbol{v}_h) + c_h(\boldsymbol{\beta}; \boldsymbol{u}_h, \boldsymbol{v}_h) + b(\boldsymbol{v}_h, p_h) - b(\boldsymbol{u}_h, q_h) = (\boldsymbol{f}, \boldsymbol{v}_h). & \text{(20b)} \end{cases}$$

Here, $\boldsymbol{u}_{0h}$ denotes an approximation of $\boldsymbol{u}_0$ belonging to $\boldsymbol{V}_h$ and the occurring forms are defined below. First note that due to (10), discretely divergence-free functions are even divergence-free pointwise; that is,

$$\boldsymbol{V}_h^{\mathrm{div}} = \{\boldsymbol{v}_h \in \boldsymbol{V}_h\colon \nabla \cdot \boldsymbol{v}_h(\boldsymbol{x}) = 0, \ \forall \boldsymbol{x} \in \Omega\}. \tag{21}$$

The solution $\boldsymbol{u}_h$ of (20) is by construction a pointwise divergence-free approximation to $\boldsymbol{u}$.

**Remark 2.3:** In the context of div-free, $\boldsymbol{H}^1$-conforming FEM, we know that the corresponding set of discretely divergence-free functions is also exactly divergence-free and thus contained in $\boldsymbol{V}^{\mathrm{div}}$. Here, $\boldsymbol{V}_h \not\subset \boldsymbol{V}$ and therefore one has to be careful: even though discretely div-free functions are div-free pointwise, we have $\boldsymbol{V}_h^{\mathrm{div}} \not\subset \boldsymbol{V}^{\mathrm{div}}$. ▲

The below appearance of certain traces of velocity facet values and normal derivatives thereof dictates that the involved velocities at least belong to $\boldsymbol{H}^{\frac{3}{2}+\varepsilon}(\mathcal{T}_h)$ for some $\varepsilon > 0$; cf. [38, Section 2.1.3]. Thus, for $\varepsilon > 0$, define the compound space

$$\boldsymbol{V}(h) = \boldsymbol{V}_h \oplus \left[\boldsymbol{V} \cap \boldsymbol{H}^{\frac{3}{2}+\varepsilon}(\mathcal{T}_h)\right]. \tag{22}$$

For the discretisation of the diffusion term, we employ the standard symmetric interior penalty (SIP) form $a_h\colon \boldsymbol{V}(h) \times \boldsymbol{V}_h \to \mathbb{R}$; cf. [38, 19]. For $\sigma > 0$, this form is given by

$$a_h(\boldsymbol{w}, \boldsymbol{v}_h) = \int_\Omega \nabla_h \boldsymbol{w} : \nabla_h \boldsymbol{v}_h \,\mathrm{d}\boldsymbol{x} - \sum_{F \in \mathcal{F}_h} \oint_F \{\!\!\{\nabla \boldsymbol{w}\}\!\!\} \boldsymbol{n}_F \cdot [\![\boldsymbol{v}_h]\!] \,\mathrm{d}\boldsymbol{s} \tag{23a}$$

$$- \sum_{F \in \mathcal{F}_h} \oint_F [\![\boldsymbol{w}]\!] \cdot \{\!\!\{\nabla \boldsymbol{v}_h\}\!\!\} \boldsymbol{n}_F \,\mathrm{d}\boldsymbol{s} + \sum_{F \in \mathcal{F}_h} \oint_F \frac{\sigma}{h_F} [\![\boldsymbol{w}]\!] \cdot [\![\boldsymbol{v}_h]\!] \,\mathrm{d}\boldsymbol{s}, \tag{23b}$$

where $[\nabla \boldsymbol{w}]_{ij} = \frac{\partial w_i}{\partial x_j}$ denotes the entries of the Jacobian. Furthermore, $h_F$ represents an appropriate length scale for the facet $F$. It is well-known that the jump penalty parameter $\sigma > 0$ has to be chosen sufficiently large such that coercivity on the discrete level is guaranteed.

**Remark 2.4:** The solution $\boldsymbol{u}_h$ to (20) is automatically normal continuous since $\boldsymbol{u}_h(t) \in \boldsymbol{V}_h \subset \boldsymbol{H}(\operatorname{div};\Omega)$ for a.e. $t \in (0,T)$. If $\boldsymbol{\tau}_F$ denotes one tangential vector to $F \in \mathcal{F}_h$, we obtain

$$\boldsymbol{v}_h = (\boldsymbol{v}_h \cdot \boldsymbol{n}_F)\boldsymbol{n}_F + (\boldsymbol{v}_h \cdot \boldsymbol{\tau}_F)\boldsymbol{\tau}_F \quad \Rightarrow \quad [\![\boldsymbol{v}_h]\!]_F = [\![(\boldsymbol{v}_h \cdot \boldsymbol{\tau}_F)]\!]\boldsymbol{\tau}_F, \quad \forall \boldsymbol{v}_h \in \boldsymbol{V}_h. \tag{24}$$

In fact, the jumps in (23) only act on tangential components of discrete velocities. ▲

**Remark 2.5:** Concerning Dirichlet (no-slip) boundary conditions, note that the normal component is prescribed in $\boldsymbol{V}_h$ to fulfil $\boldsymbol{u}_h \cdot \boldsymbol{n}\big|_{\partial\Omega} = 0$—the no-penetration condition. Thus, the boundary facet contribution for $\mathcal{F}_h^{\partial} \subset \mathcal{F}_h$ in (23) only acts on tangential components. This weak imposition of the no-slip condition is therefore consistent with the limiting case of $\nu \to 0$. Even more, for high Reynolds numbers, imposing Dirichlet boundary conditions by means of a Nitsche penalty method can be considered as an implicit wall model [23]. To this end, the parameter $\sigma$ can also be designed using Spalding's law of the wall [4]. Thus, in certain situations, it may be advantageous to assign a different value of $\sigma$ to no-slip facets $F \in \mathcal{F}_h^{\partial}$ than to interior facets $F \in \mathcal{F}_h^{i}$. ▲

In conjunction with the viscous term $a_h$, the following expressions are used:

$$\forall \boldsymbol{w} \in \boldsymbol{V}(h): \quad |\!|\!|\boldsymbol{w}|\!|\!|_e^2 = \|\nabla_h \boldsymbol{w}\|_{\boldsymbol{L}^2}^2 + \sum_{F \in \mathcal{F}_h} \frac{\sigma}{h_F} \|[\![\boldsymbol{w}]\!]\|_{\boldsymbol{L}^2(F)}^2 \tag{25a}$$

$$|\!|\!|\boldsymbol{w}|\!|\!|_{e,\sharp}^2 = |\!|\!|\boldsymbol{w}|\!|\!|_e^2 + \sum_{K \in \mathcal{T}_h} h_K \|\nabla \boldsymbol{w} \cdot \boldsymbol{n}_K\|_{\boldsymbol{L}^2(\partial K)}^2 \tag{25b}$$

Here, $|\!|\!|\cdot|\!|\!|_e$ denotes a discrete energy norm and the index $\sharp$ indicates a stronger norm.

The pressure-velocity coupling $b\colon \boldsymbol{V}(h) \times Q \to \mathbb{R}$ remains unchanged:

$$b(\boldsymbol{w},q) = -\int_\Omega q(\nabla \cdot \boldsymbol{w})\,\mathrm{d}\boldsymbol{x} \tag{26}$$

For the (linearised) inertia term, we choose a convection term $c_h\colon \boldsymbol{V}(h) \times \boldsymbol{V}_h \to \mathbb{R}$ with optional upwinding controlled by the parameter $\gamma \geqslant 0$:

$$c_h(\boldsymbol{\beta};\boldsymbol{w},\boldsymbol{v}_h) = \int_\Omega (\boldsymbol{\beta} \cdot \nabla_h)\boldsymbol{w} \cdot \boldsymbol{v}_h \,\mathrm{d}\boldsymbol{x} - \sum_{F \in \mathcal{F}_h^i} \oint_F (\boldsymbol{\beta} \cdot \boldsymbol{n}_F)[\![\boldsymbol{w}]\!] \cdot \{\!\{\boldsymbol{v}_h\}\!\}\,\mathrm{d}\boldsymbol{s} \tag{27a}$$

$$+ \sum_{F \in \mathcal{F}_h^i} \oint_F \frac{\gamma}{2} |\boldsymbol{\beta} \cdot \boldsymbol{n}_F| [\![\boldsymbol{w}]\!] \cdot [\![\boldsymbol{v}_h]\!]\,\mathrm{d}\boldsymbol{s} \tag{27b}$$

Note that due to the strong imposition of the no-penetration condition $\boldsymbol{\beta} \cdot \boldsymbol{n}\big|_{\partial\Omega} = 0$ in $\boldsymbol{V}_h$, the convective form $c_h$ does not contain any surface integrals over boundary facets. Therefore, as long as there is no in- or outflow across $\partial\Omega$, the weak imposition of tangential boundary conditions is handled exclusively by $a_h$. Moreover, a modification of $c_h$ which ensures skew-symmetry is redundant since $\nabla \cdot \boldsymbol{\beta} = 0$ has been assumed.

**Remark 2.6:** Let us comment on the upwind stabilisation in (27). This kind of convection stabilisation, as opposed to dissipative viscous stabilisation, does not explicitly include gradient-based terms. However, as can be seen in [19, Section 2.3] for the stationary transport problem, upwind stabilisation gives additional control over the streamline derivative—much like classical SUPG stabilisation. Unfortunately, we are not (yet) aware of how to show (if at all) an analogous result can be obtained for incompressible flow problems either in the full dG-FEM setting, or in the $\boldsymbol{H}(\mathrm{div})$-FEM case. ▲

In conjunction with the convection term $c_h$, the following expressions are used:

$$\forall \boldsymbol{w} \in \boldsymbol{V}(h): \quad |||\boldsymbol{w}|||_{\boldsymbol{\beta}}^2 = \|\boldsymbol{w}\|_{\boldsymbol{L}^2}^2 + |\boldsymbol{w}|_{\boldsymbol{\beta},\mathrm{upw}}^2, \quad |\boldsymbol{w}|_{\boldsymbol{\beta},\mathrm{upw}}^2 = \sum_{F\in\mathcal{F}_h^i} \oint_F \frac{\gamma}{2} |\boldsymbol{\beta}\cdot\boldsymbol{n}_F| |[\![\boldsymbol{w}]\!]|^2 \,\mathrm{d}\boldsymbol{s} \tag{28a}$$

$$|||\boldsymbol{w}|||_{\boldsymbol{\beta},\sharp}^2 = |||\boldsymbol{w}|||_{\boldsymbol{\beta}}^2 + \|\boldsymbol{\beta}\|_\infty^2 \sum_{K\in\mathcal{T}_h} h_K^{-2} \|\boldsymbol{w}\|_{\boldsymbol{L}^2(K)}^2 + \|\boldsymbol{\beta}\|_\infty \sum_{K\in\mathcal{T}_h} \|\boldsymbol{w}\|_{\boldsymbol{L}^2(\partial K)}^2 \tag{28b}$$

Here, $|||\cdot|||_{\boldsymbol{\beta}}$ measures terms due to convection and again, the index $\sharp$ indicates a stronger norm since non-zero terms are added. Moreover, $|\cdot|_{\boldsymbol{\beta},\mathrm{upw}}$ denotes the upwind seminorm, which represents additional control over $\boldsymbol{\beta}$-scaled velocity jumps.

**Remark 2.7:** Concerning the appearing forms in (20), one could alternatively take any other inf-sup stable dG-FEM formulation (cf., for example, [38, 8]), and adapt it by neglecting all terms with $\boldsymbol{v}_h\cdot\boldsymbol{n}\big|_{\partial\Omega}$ and $[\![\boldsymbol{v}_h]\!]\cdot\boldsymbol{n}_F$ for $F\in\mathcal{F}_h^i$. In our setting, these terms vanish automatically due to the strong imposition of the no-penetration condition and the continuity of normal components of functions in $\boldsymbol{H}(\mathrm{div};\Omega)$. ▲

**Remark 2.8:** Until now we only considered no-slip conditions. However, the definitions of $a_h$ and $c_h$ can be extended directly to the case of weakly imposed periodic boundary conditions (BCs). To this end, instead of perceiving the set of periodic facets as boundary facets, it is very natural to treat them analogously to interior facets. In order to keep the $\boldsymbol{H}(\mathrm{div})$-conformity intact, only the normal continuity across the periodic boundary has to be ensured in a strong sense. Thus, the subsequent analysis holds true *verbatim* for problems involving periodic BCs. We refer to [47] where periodic BCs for scalar diffusion-reaction problem are considered in the dG-FEM context. ▲

In the following, the Galerkin orthogonality property of (20) is stated.

**Corollary 2.9 (*Galerkin orthogonality*)**

*Let $\boldsymbol{u}_h \in \boldsymbol{V}_h^T$ solve (20), and assume that the solution $\boldsymbol{u}\in\boldsymbol{V}^T$ of (3) satisfies the minimum regularity $\boldsymbol{u}\in L^2\left(0,T;\boldsymbol{H}^{\frac{3}{2}+\varepsilon}(\mathcal{T}_h)\right)$ for $\varepsilon>0$. Then, for all $\boldsymbol{v}_h\in\boldsymbol{V}_h^{\mathrm{div}}$:*

$$(\partial_t[\boldsymbol{u}-\boldsymbol{u}_h],\boldsymbol{v}_h) + \nu a_h(\boldsymbol{u}-\boldsymbol{u}_h,\boldsymbol{v}_h) + c_h(\boldsymbol{\beta};\boldsymbol{u}-\boldsymbol{u}_h,\boldsymbol{v}_h) = 0, \quad \text{for a.e. } t\in(0,T) \tag{29}$$

**Proof:** The most important ingredient is the consistency of both SIP formulation of the viscous term and upwind formulation of the convective term. We will not prove this here, but instead refer to [38, 19]. Having the consistency property in mind, subtracting (20) from (3) and using arbitrary $(\boldsymbol{v}_h,q_h)\in\boldsymbol{V}_h\times Q_h$ as test functions leads to

$$(\partial_t[\boldsymbol{u}-\boldsymbol{u}_h],\boldsymbol{v}_h) + \nu a_h(\boldsymbol{u}-\boldsymbol{u}_h,\boldsymbol{v}_h) + c_h(\boldsymbol{\beta};\boldsymbol{u}-\boldsymbol{u}_h,\boldsymbol{v}_h) \tag{30a}$$

$$+ b(\boldsymbol{v}_h,p-p_h) - b(\boldsymbol{u}-\boldsymbol{u}_h,q_h) = 0. \tag{30b}$$

Restricting the velocity test functions to discretely divergence-free ones, that is $\boldsymbol{v}_h \in \boldsymbol{V}_h^{\text{div}}$, we can use $\nabla \cdot \boldsymbol{v}_h = 0$ to remove the first mixed term; see (21). The second mixed term vanishes by construction because both $\boldsymbol{u}$ and $\boldsymbol{u}_h$ are in $\boldsymbol{V}_h^{\text{div}}$ and $Q_h \subset Q$. ∎

In order to obtain stability estimates for (20), the following results have to be established.

**Lemma 2.10** (***Discrete coercivity of*** $a_h$ ***and*** $c_h$)

*Assume that $\sigma > 0$ is sufficiently large. Then, the SIP bilinear form $a_h$ is coercive on $\boldsymbol{V}_h$ w.r.t. the energy norm $|||\cdot|||_e$. Moreover, the convective form $c_h$ is coercive on $\boldsymbol{V}_h$ w.r.t. the upwind seminorm $|\cdot|_{\boldsymbol{\beta},\text{upw}}$. That is, there exists $C_\sigma > 0$, independent of $h$, such that,*

$$\forall \boldsymbol{v}_h \in \boldsymbol{V}_h: \quad a_h(\boldsymbol{v}_h, \boldsymbol{v}_h) \geqslant C_\sigma |||\boldsymbol{v}_h|||_e^2 \quad \text{and} \quad c_h(\boldsymbol{\beta}; \boldsymbol{v}_h, \boldsymbol{v}_h) = |\boldsymbol{v}_h|_{\boldsymbol{\beta},\text{upw}}^2. \tag{31}$$

**Proof:** Cf., for example, [38, Lemma 6.6] or [19, Section 6.1.2.1] for the discrete coercivity of $a_h$. The non-negativity of $c_h$ is shown, for example, in [19, Lemma 6.39]. ∎

As already said, we want to decouple velocity and pressure and treat them both separately. For this separation to work it is essential to work in $\boldsymbol{V}_h^{\text{div}}$ which, due to (11), is non-trivial.

**Lemma 2.11** (***Well-posedness and velocity energy estimate***)

*Let $\boldsymbol{f} \in L^1(\boldsymbol{L}^2)$ and $\boldsymbol{u}_{0h} \in \boldsymbol{L}^2$. Then, there exists a solution $\boldsymbol{u}_h \in \boldsymbol{V}_h^T$ to (20) with*

$$\frac{1}{2} \|\boldsymbol{u}_h\|_{L^\infty(\boldsymbol{L}^2)} + \int_0^T \left[ \nu C_\sigma |||\boldsymbol{u}_h|||_e^2 + |\boldsymbol{u}_h|_{\boldsymbol{\beta},\text{upw}}^2 \right] \mathrm{d}\tau \leqslant \|\boldsymbol{u}_{0h}\|_{\boldsymbol{L}^2}^2 + \frac{3}{2} \|\boldsymbol{f}\|_{L^1(\boldsymbol{L}^2)}^2. \tag{32}$$

*Provided $\boldsymbol{f}$ is even Lipschitz in time, the solution $\boldsymbol{u}_h$ is unique.*

**Proof:** Testing (20) with $(\boldsymbol{u}_h(t), 0) \in \boldsymbol{V}_h^{\text{div}} \times Q_h$, using the discrete coercivity properties in Lemma 2.10 on the left-hand side and Cauchy–Schwarz on the right-hand side leads to

$$\frac{1}{2}\frac{\mathrm{d}}{\mathrm{d}t} \|\boldsymbol{u}_h(t)\|_{\boldsymbol{L}^2}^2 + \nu C_\sigma |||\boldsymbol{u}_h(t)|||_e^2 + |\boldsymbol{u}_h(t)|_{\boldsymbol{\beta},\text{upw}}^2 \leqslant \|\boldsymbol{f}(t)\|_{\boldsymbol{L}^2} \|\boldsymbol{u}_h(t)\|_{\boldsymbol{L}^2}, \tag{33}$$

since $(\partial_t \boldsymbol{u}_h, \boldsymbol{u}_h) = \frac{1}{2}\frac{\mathrm{d}}{\mathrm{d}t} \|\boldsymbol{u}_h\|_{\boldsymbol{L}^2}^2$. Now, we directly follow [17], where the estimate

$$\|\boldsymbol{u}_h(t)\|_{\boldsymbol{L}^2} \leqslant \|\boldsymbol{u}_{0h}\|_{\boldsymbol{L}^2} + \|\boldsymbol{f}\|_{L^1(\boldsymbol{L}^2)} \tag{34}$$

plays a key role. Inserting (34) in (33), applying Young's inequality and integrating over $(0,T)$ shows the estimate. Applying the theorem of Carathéodory (cf. [31, Theorem A.50]) ensures existence and, if $\boldsymbol{f}$ is additionally Lipschitz in time, uniqueness. ∎

**Remark 2.12:** Alternative results can be obtained by estimating the right-hand side of (33) differently. Firstly, one could use Young's and Poincare's inequality with the goal of bounding $\|\boldsymbol{u}_h\|_{\boldsymbol{L}^2}$ in an appropriately scaled energy norm $|||\boldsymbol{u}_h|||_e$. This inevitably leads to a $\nu^{-1}$ factor. Secondly, also after Young's inequality, the Gronwall lemma could be applied. A factor $\exp(T)$ appears on the right-hand side. ▲

## 3 Pressure- and *Re*-semi-robust analysis for the linearised problem

From now on, $C > 0$ denotes a generic constant independent of $h$ and $\nu$. The first subsection deals with optimal-order, pressure- and *Re*-semi-robust estimates for the discrete, stationary Stokes projection in the $\boldsymbol{H}(\mathrm{div})$ context. This is an important step, since, for the analysis of the Oseen problem in the second subsection, the stationary Stokes projection is used for the error splitting. In this way, the approximation properties of the projection operator in the Oseen problem can be derived from error estimates for the stationary Stokes problem.

### 3.1 Stationary Stokes projection

In this section we basically consider the stationary Stokes problem. With a sufficiently smooth forcing term $\boldsymbol{g}$, the well-known continuous weak formulation reads

$$\begin{cases} ?\,(\boldsymbol{u}_s, p_s) \in \boldsymbol{V} \times Q \text{ s.t.}, \forall (\boldsymbol{v}, q) \in \boldsymbol{V} \times Q, & \text{(35a)} \\ \nu a(\boldsymbol{u}_s, \boldsymbol{v}) + b(\boldsymbol{v}, p_s) - b(\boldsymbol{u}_s, q) = (\boldsymbol{g}, \boldsymbol{v}). & \text{(35b)} \end{cases}$$

In order to obtain optimal $\boldsymbol{L}^2$-estimates for the velocity, we make the following assumption which is called 'elliptic regularity', 'Cattabriga's regularity' or 'smoothing property'.

**Assumption A7:** Assume that $\Omega$ is either a convex polygon for $d = 2$ or of class $\mathcal{C}^{1,1}$ for $d \in \{2,3\}$. Then, for all $\boldsymbol{g} \in \boldsymbol{L}^2$, the solution $(\boldsymbol{u}_s, p_s) \in \boldsymbol{V} \times Q$ of (35) additionally fulfils the regularity property $(\boldsymbol{u}_s, p_s) \in \boldsymbol{H}^2 \times H^1$ and the *a priori* estimate $\sqrt{\nu}\,\|\boldsymbol{u}_s\|_{\boldsymbol{H}^2} + \|p_s\|_{H^1} \leqslant C\,\|\boldsymbol{g}\|_{\boldsymbol{L}^2}$; cf. [9, Theorem IV.5.8].

Extending [21,31] to the $\boldsymbol{H}(\mathrm{div})$-conforming case, we give the following definition. Note that the definition is stated directly in $\boldsymbol{V}_h^{\mathrm{div}}$ because this suffices for our considerations.

**Definition 3.1 (*Stationary Stokes projection*)**

Let $\boldsymbol{w} \in \boldsymbol{H}^{\frac{3}{2}+\varepsilon}(\mathcal{T}_h)$ for $\varepsilon > 0$ fulfil $\nabla \cdot \boldsymbol{w} = 0$ pointwise. Then, we define the stationary Stokes projection $\boldsymbol{\pi}_s \boldsymbol{w} \in \boldsymbol{V}_h^{\mathrm{div}}$ of $\boldsymbol{w}$ to be the unique FE solution to the problem

$$a_h(\boldsymbol{\pi}_s \boldsymbol{w}, \boldsymbol{v}_h) = a_h(\boldsymbol{w}, \boldsymbol{v}_h), \quad \forall \boldsymbol{v}_h \in \boldsymbol{V}_h^{\mathrm{div}}. \tag{36}$$

Analogously to divergence-free $\boldsymbol{H}^1$-conforming FEM, the stationary Stokes projection coincides with a vector-valued elliptic (or Ritz) projection on $\boldsymbol{V}_h^{\mathrm{div}}$; cf. [42]. Very conveniently, the approximation properties of the projection operator $\boldsymbol{\pi}_s$ can thus be derived from error estimates for the stationary Stokes problem.

Additionally to coercivity of $a_h$ (Lemma 2.10) the following continuity result is required.

**Lemma 3.2 (*Boundedness of diffusion term* $a_h$)**

*There exists a* $M_{\mathrm{dif}} > 0$, *independent of* $h$, *such that*

$$\forall (\boldsymbol{w}, \boldsymbol{v}_h) \in \boldsymbol{V}(h) \times \boldsymbol{V}_h: \quad a_h(\boldsymbol{w}, \boldsymbol{v}_h) \leqslant M_{\mathrm{dif}} |||\boldsymbol{w}|||_{e,\sharp} |||\boldsymbol{v}_h|||_e. \tag{37}$$

**Proof:** Cf., for example, [19, Section 4.2.3] for a scalar-valued analogue. Up to positive factors independent of $h$, the claim follows from a componentwise application. ■

**THEOREM 3.3** (*Stokes projection error estimate*)

*Let $\boldsymbol{\pi}_s \boldsymbol{w}$ be the Stokes projection of $\boldsymbol{w}$ and assume elliptic regularity A7. Then,*

$$\|\boldsymbol{w}-\boldsymbol{\pi}_s\boldsymbol{w}\|_{\boldsymbol{L}^2} + h|||\boldsymbol{w}-\boldsymbol{\pi}_s\boldsymbol{w}|||_{e,\sharp} \leqslant C_s h \inf_{\boldsymbol{w}_h \in \boldsymbol{V}_h^{\mathrm{div}}} |||\boldsymbol{w}-\boldsymbol{w}_h|||_{e,\sharp}. \tag{38}$$

**PROOF:** (i) We begin with a pressure-robust estimate in the $|||\cdot|||_{e,\sharp}$-norm. Let $\boldsymbol{w}_h \in \boldsymbol{V}_h^{\mathrm{div}}$ be arbitrary and define $\boldsymbol{w}_h^0 = \boldsymbol{\pi}_s\boldsymbol{w} - \boldsymbol{w}_h \in \boldsymbol{V}_h^{\mathrm{div}}$. The Galerkin orthogonality being inherent in Definition 3.1 yields

$$a_h\big(\boldsymbol{w}_h^0, \boldsymbol{v}_h\big) = a_h(\boldsymbol{\pi}_s\boldsymbol{w}-\boldsymbol{w}_h, \boldsymbol{v}_h) = a_h(\boldsymbol{w}-\boldsymbol{w}_h, \boldsymbol{v}_h), \quad \forall \boldsymbol{v}_h \in \boldsymbol{V}_h^{\mathrm{div}}. \tag{39}$$

Now, choose $\boldsymbol{v}_h = \boldsymbol{w}_h^0 \in \boldsymbol{V}_h^{\mathrm{div}}$, use discrete coercivity of $a_h$ (Lemma 2.10) on the left-hand side and Lemma 3.2 and Young's inequality ($\varepsilon > 0$) on the right-hand side to obtain

$$C_\sigma \big|\big|\big|\boldsymbol{w}_h^0\big|\big|\big|_e^2 \leqslant M_{\mathrm{dif}} |||\boldsymbol{w}-\boldsymbol{w}_h|||_{e,\sharp} \big|\big|\big|\boldsymbol{w}_h^0\big|\big|\big|_e \leqslant \frac{1}{2\varepsilon} M_{\mathrm{dif}} |||\boldsymbol{w}-\boldsymbol{w}_h|||_{e,\sharp}^2 + \frac{\varepsilon}{2} M_{\mathrm{dif}} \big|\big|\big|\boldsymbol{w}_h^0\big|\big|\big|_e^2. \tag{40}$$

Choosing $\varepsilon = C_\sigma M_{\mathrm{dif}}^{-1}$, multiplication by 2 and reordering leads to

$$C_\sigma |||\boldsymbol{\pi}_s\boldsymbol{w}-\boldsymbol{w}_h|||_e^2 \leqslant \frac{M_{\mathrm{dif}}^2}{C_\sigma} |||\boldsymbol{u}-\boldsymbol{w}_h|||_{e,\sharp}^2. \tag{41}$$

Because $\boldsymbol{w}_h \in \boldsymbol{V}_h^{\mathrm{div}}$ is arbitrary, the triangle inequality and $|||\cdot|||_e \leqslant |||\cdot|||_{e,\sharp}$ yield

$$|||\boldsymbol{w}-\boldsymbol{\pi}_s\boldsymbol{w}|||_e^2 \leqslant 2\Big[|||\boldsymbol{w}-\boldsymbol{w}_h|||_e^2 + |||\boldsymbol{\pi}_s\boldsymbol{w}-\boldsymbol{w}_h|||_e^2\Big] \leqslant C \inf_{\boldsymbol{w}_h \in \boldsymbol{V}_h^{\mathrm{div}}} |||\boldsymbol{w}-\boldsymbol{w}_h|||_{e,\sharp}^2. \tag{42}$$

The final step is to acknowledge that the $|||\cdot|||_e$- and $|||\cdot|||_{e,\sharp}$-norm are uniformly equivalent on $\boldsymbol{V}_h$, that is, there exists a $C > 0$ such that $C|||\cdot|||_{e,\sharp} \leqslant |||\cdot|||_e$; cf. [19]. Thus,

$$|||\boldsymbol{w}-\boldsymbol{\pi}_s\boldsymbol{w}|||_{e,\sharp}^2 \leqslant C \inf_{\boldsymbol{w}_h \in \boldsymbol{V}_h^{\mathrm{div}}} |||\boldsymbol{w}-\boldsymbol{w}_h|||_{e,\sharp}^2. \tag{43}$$

Note that we did not include the viscosity in the definition (36) of the Stokes projection. Hence, the constant $C$ in the last estimate is not corrupted by negative powers of $\nu$.

(ii) Secondly, in order to obtain an estimate also for the kinetic energy, the elliptic regularity condition is essential. For brevity, we omit to give full details at this point. However, a careful inspection of [19, Theorem 6.25] reveals that in our case, it is possible to obtain the following pressure-robust estimate:

$$\|\boldsymbol{w}-\boldsymbol{\pi}_s\boldsymbol{w}\|_{\boldsymbol{L}^2}^2 \leqslant Ch^2 |||\boldsymbol{w}-\boldsymbol{\pi}_s\boldsymbol{w}|||_{e,\sharp}^2 \leqslant Ch^2 \inf_{\boldsymbol{w}_h \in \boldsymbol{V}_h^{\mathrm{div}}} |||\boldsymbol{w}-\boldsymbol{w}_h|||_{e,\sharp}^2 \tag{44}$$

At this point, it is important that we are dealing with a divergence-free method; cf. [22, Theorem 6.4.2] for an $\boldsymbol{H}^1$-conforming exposition. Otherwise, the $\boldsymbol{L}^2$-estimate would not readily be pressure-robust. ∎

### 3.2 Velocity error estimates for the Oseen problem

In the previous subsection, we provided estimates for the viscous Stokes term. Now, we present an estimate for the Oseen convection term which allows for *Re*-semi-robust estimates.

**Lemma 3.4 (*Boundedness of convection term $c_h$*)**

*Let $M_{\mathrm{cnv}} = \max\left\{1, C_{\mathrm{inv}}^2, \gamma^{-1}\right\}$ and $\varepsilon_1, \varepsilon_2 > 0$. Then, for all $(\boldsymbol{w}, \boldsymbol{v}_h) \in \boldsymbol{V}(h) \times \boldsymbol{V}_h$,*

$$|c_h(\boldsymbol{\beta}; \boldsymbol{w}, \boldsymbol{v}_h)| \leqslant M_{\mathrm{cnv}} \left( \frac{1}{2\varepsilon_1} + \frac{1}{\varepsilon_2} \right) |||\boldsymbol{w}|||_{\boldsymbol{\beta},\sharp}^2 + \frac{\varepsilon_1}{2} \|\boldsymbol{v}_h\|_{\boldsymbol{L}^2}^2 + \varepsilon_2 |\boldsymbol{v}_h|_{\boldsymbol{\beta},\mathrm{upw}}^2. \tag{45}$$

**Proof:** We begin with the following integration by parts variant:

$$\sum_{F \in \mathcal{F}_h} \oint_F (\boldsymbol{\beta} \cdot \boldsymbol{n}_F) \Big[ [\![\boldsymbol{w}]\!] \cdot \{\!\{\boldsymbol{v}_h\}\!\} + [\![\boldsymbol{v}_h]\!] \cdot \{\!\{\boldsymbol{w}\}\!\} \Big] \,\mathrm{d}\boldsymbol{s} \tag{46a}$$

$$= \sum_{K \in \mathcal{T}_h} \oint_{\partial K} (\boldsymbol{\beta} \cdot \boldsymbol{n}_K) \boldsymbol{w} \cdot \boldsymbol{v}_h \,\mathrm{d}\boldsymbol{s} = \sum_{K \in \mathcal{T}_h} \int_K \Big[ (\boldsymbol{\beta} \cdot \nabla) \boldsymbol{w} \cdot \boldsymbol{v}_h + (\boldsymbol{\beta} \cdot \nabla) \boldsymbol{v}_h \cdot \boldsymbol{w} \Big] \,\mathrm{d}\boldsymbol{x} \tag{46b}$$

The first equality is a result of $[\![\boldsymbol{w} \cdot \boldsymbol{v}_h]\!] = [\![\boldsymbol{w}]\!] \cdot \{\!\{\boldsymbol{v}_h\}\!\} + [\![\boldsymbol{v}_h]\!] \cdot \{\!\{\boldsymbol{w}\}\!\}$ and the second equality can be shown using elementwise integration by parts with $\nabla \cdot \boldsymbol{\beta} = 0$ and the product rule for the gradient. Inserting this into the convective form (27) leads to

$$c_h(\boldsymbol{\beta}; \boldsymbol{w}, \boldsymbol{v}_h) = -\sum_{K \in \mathcal{T}_h} \int_K (\boldsymbol{\beta} \cdot \nabla) \boldsymbol{v}_h \cdot \boldsymbol{w} \,\mathrm{d}\boldsymbol{x} + \sum_{F \in \mathcal{F}_h^i} \oint_F (\boldsymbol{\beta} \cdot \boldsymbol{n}_F) [\![\boldsymbol{v}_h]\!] \cdot \{\!\{\boldsymbol{w}\}\!\} \,\mathrm{d}\boldsymbol{s} \tag{47a}$$

$$+ \sum_{F \in \mathcal{F}_h^i} \oint_F \frac{\gamma}{2} |\boldsymbol{\beta} \cdot \boldsymbol{n}_F| [\![\boldsymbol{v}_h]\!] \cdot [\![\boldsymbol{w}]\!] \,\mathrm{d}\boldsymbol{s} = \mathfrak{T}_1 + \mathfrak{T}_2 + \mathfrak{T}_3. \tag{47b}$$

For the first term, using the generalised Hölder inequality, the local inverse inequality, Cauchy–Schwarz and Young's inequality with $\varepsilon_1 > 0$, we infer

$$|\mathfrak{T}_1| \leqslant \sum_{K \in \mathcal{T}_h} \|\boldsymbol{\beta}\|_{\boldsymbol{L}^\infty(K)} \|\nabla \boldsymbol{v}_h\|_{\boldsymbol{L}^2(K)} \|\boldsymbol{w}\|_{\boldsymbol{L}^2(K)} \tag{48a}$$

$$\leqslant \sum_{K \in \mathcal{T}_h} C_{\mathrm{inv}} h_K^{-1} \|\boldsymbol{\beta}\|_{\boldsymbol{L}^\infty(K)} \|\boldsymbol{v}_h\|_{\boldsymbol{L}^2(K)} \|\boldsymbol{w}\|_{\boldsymbol{L}^2(K)} \tag{48b}$$

$$\leqslant \left( \sum_{K \in \mathcal{T}_h} \|\boldsymbol{v}_h\|_{\boldsymbol{L}^2(K)}^2 \right)^{1/2} C_{\mathrm{inv}} \left( \sum_{K \in \mathcal{T}_h} h_K^{-2} \|\boldsymbol{\beta}\|_{\boldsymbol{L}^\infty(K)}^2 \|\boldsymbol{w}\|_{\boldsymbol{L}^2(K)}^2 \right)^{1/2} \tag{48c}$$

$$\leqslant \|\boldsymbol{v}_h\|_{\boldsymbol{L}^2} C_{\mathrm{inv}} |||\boldsymbol{w}|||_{\boldsymbol{\beta},\sharp} \leqslant \frac{\varepsilon_1}{2} \|\boldsymbol{v}_h\|_{\boldsymbol{L}^2}^2 + \frac{1}{2\varepsilon_1} C_{\mathrm{inv}}^2 |||\boldsymbol{w}|||_{\boldsymbol{\beta},\sharp}^2. \tag{48d}$$

For the first facet term, we use the definition of the absolute value, Cauchy–Schwarz and Young's inequality ($\varepsilon_2 > 0$) to obtain

$$\mathfrak{T}_2 \leqslant \sum_{F\in\mathcal{F}_h^i} \oint_F \left[ \left(\frac{\gamma}{2}|\boldsymbol{\beta}\cdot\boldsymbol{n}_F|\right)^{1/2} [\![\boldsymbol{v}_h]\!] \cdot \left(\frac{2}{\gamma}|\boldsymbol{\beta}\cdot\boldsymbol{n}_F|\right)^{1/2} \{\!\!\{\boldsymbol{w}\}\!\!\} \right] \mathrm{d}\boldsymbol{s} \tag{49a}$$

$$\leqslant |\boldsymbol{v}_h|_{\boldsymbol{\beta},\mathrm{upw}} \left( \sum_{F\in\mathcal{F}_h^i} \oint_F \frac{2}{\gamma}|\boldsymbol{\beta}\cdot\boldsymbol{n}_F| \left|\{\!\!\{\boldsymbol{w}\}\!\!\}\right|^2 \mathrm{d}\boldsymbol{s} \right)^{1/2} \tag{49b}$$

$$\leqslant |\boldsymbol{v}_h|_{\boldsymbol{\beta},\mathrm{upw}} \left( \frac{1}{\gamma} \sum_{K\in\mathcal{T}_h} \|\boldsymbol{\beta}\|_{\boldsymbol{L}^\infty(K)} \|\boldsymbol{w}\|^2_{\boldsymbol{L}^2(\partial K)} \right)^{1/2} \tag{49c}$$

$$\leqslant |\boldsymbol{v}_h|_{\boldsymbol{\beta},\mathrm{upw}} \gamma^{-1/2} |\!|\!|\boldsymbol{w}|\!|\!|_{\boldsymbol{\beta},\sharp} \leqslant \frac{\varepsilon_2}{2} |\boldsymbol{v}_h|^2_{\boldsymbol{\beta},\mathrm{upw}} + \frac{1}{2\varepsilon_2}\gamma^{-1} |\!|\!|\boldsymbol{w}|\!|\!|^2_{\boldsymbol{\beta},\sharp}, \tag{49d}$$

where the third step uses the bound $(a+b)^2 \leqslant 2(a^2+b^2)$ for $a,b\in\mathbb{R}$ which yields

$$\sum_{F\in\mathcal{F}_h^i} \oint_F 2\left|\{\!\!\{\boldsymbol{w}\}\!\!\}\right|^2 \mathrm{d}\boldsymbol{s} \leqslant \sum_{F\in\mathcal{F}_h^i} \left[ \|\boldsymbol{w}^+\|^2_{\boldsymbol{L}^2(F)} + \|\boldsymbol{w}^-\|^2_{\boldsymbol{L}^2(F)} \right] \leqslant \sum_{K\in\mathcal{T}_h} \|\boldsymbol{w}\|^2_{\boldsymbol{L}^2(\partial K)}. \tag{50}$$

With Young's inequality ($\varepsilon_2 > 0$), the upwind term can be bounded trivially by

$$\mathfrak{T}_3 \leqslant |\boldsymbol{v}_h|_{\boldsymbol{\beta},\mathrm{upw}} |\boldsymbol{w}|_{\boldsymbol{\beta},\mathrm{upw}} \leqslant |\boldsymbol{v}_h|_{\boldsymbol{\beta},\mathrm{upw}} |\!|\!|\boldsymbol{w}|\!|\!|_{\boldsymbol{\beta},\sharp} \leqslant \frac{\varepsilon_2}{2}|\boldsymbol{v}_h|^2_{\boldsymbol{\beta},\mathrm{upw}} + \frac{1}{2\varepsilon_2} |\!|\!|\boldsymbol{w}|\!|\!|_{\boldsymbol{\beta},\sharp}. \tag{51}$$

Collecting the above estimates concludes the proof. ■

**Remark 3.5:** Note that in [27], $\boldsymbol{H}(\mathrm{div})$-FEM for the incompressible Euler equations ($\nu = 0$) have been considered. There, one can find an estimate for the corresponding difference of convective terms in the nonlinear case $\boldsymbol{\beta} = \boldsymbol{u}$. However, it has to be mentioned that the analysis in [27] heavily relies on the regularity assumption $\hat{\boldsymbol{u}} \in \boldsymbol{W}^{1,\infty}([0,T]\times\Omega)$ for the solution $\hat{\boldsymbol{u}}$ of the incompressible Euler problem; see [27, Theorem 2.2]. This assumption is very restrictive in the case of $\nu = 0$ as there exists no inherent smoothing mechanism from the incompressible Euler operator in the crosswind direction(s). ▲

Now, we can use the Stokes projection to introduce an error splitting:

$$\boldsymbol{u} - \boldsymbol{u}_h = [\boldsymbol{u} - \boldsymbol{\pi}_s\boldsymbol{u}] - [\boldsymbol{u}_h - \boldsymbol{\pi}_s\boldsymbol{u}] = \boldsymbol{\eta} - \boldsymbol{e}_h \tag{52}$$

**Theorem 3.6** (***Velocity discretisation error estimate***)

*Let $\boldsymbol{u} \in \boldsymbol{V}^T$ solve (3) and $\boldsymbol{u}_h \in \boldsymbol{V}_h^T$ solve (20). If additionally $\boldsymbol{u} \in L^2\left(\boldsymbol{H}^{\frac{3}{2}+\varepsilon}(\mathcal{T}_h)\right)$ for $\varepsilon > 0$, $\boldsymbol{\beta} \in L^\infty(\boldsymbol{L}^\infty)$ and $\boldsymbol{u}_h(0) = \boldsymbol{\pi}_s\boldsymbol{u}_0$, for arbitrary $\alpha > 0$, we obtain:*

$$\|\boldsymbol{e}_h\|^2_{L^\infty(\boldsymbol{L}^2)} + \int_0^T \left[ \nu C_\sigma |\!|\!|\boldsymbol{e}_h(\tau)|\!|\!|^2_e + |\boldsymbol{e}_h(\tau)|^2_{\boldsymbol{\beta},\mathrm{upw}} \right] \mathrm{d}\tau \tag{53a}$$

$$\leqslant C\alpha^{-1} e^{\alpha T} \int_0^T \left[ \|\partial_t \boldsymbol{\eta}(\tau)\|^2_{\boldsymbol{L}^2} + |\!|\!|\boldsymbol{\eta}(\tau)|\!|\!|^2_{\boldsymbol{\beta},\sharp} \right] \mathrm{d}\tau \tag{53b}$$

**Proof:** Galerkin orthogonality (Corollary 2.9) with $\boldsymbol{v}_h = \boldsymbol{e}_h(t) \in \boldsymbol{V}_h^{\mathrm{div}}$ and (52) yields

$$(\partial_t \boldsymbol{e}_h, \boldsymbol{e}_h) + \nu a_h(\boldsymbol{e}_h, \boldsymbol{e}_h) + c_h(\boldsymbol{\beta}; \boldsymbol{e}_h, \boldsymbol{e}_h) = (\partial_t \boldsymbol{\eta}, \boldsymbol{e}_h) + \nu a_h(\boldsymbol{\eta}, \boldsymbol{e}_h) + c_h(\boldsymbol{\beta}; \boldsymbol{\eta}, \boldsymbol{e}_h). \tag{54}$$

We use $(\partial_t \boldsymbol{e}_h, \boldsymbol{e}_h) = \frac{1}{2}\frac{\mathrm{d}}{\mathrm{d}t}\|\boldsymbol{e}_h\|_{\boldsymbol{L}^2}^2$ and discrete coercivity of $a_h$ and $c_h$ (Lemma 2.10) on the left-hand side. On the right-hand side, we apply Cauchy–Schwarz plus Young ($\varepsilon_3 > 0$) and the definition for the stationary Stokes projection (Definition 3.1). Then, boundedness of the convective term (Lemma 3.4) leads to

$$\frac{1}{2}\frac{\mathrm{d}}{\mathrm{d}t}\|\boldsymbol{e}_h\|_{\boldsymbol{L}^2}^2 + \nu C_\sigma |\!|\!| \boldsymbol{e}_h |\!|\!|_e^2 + |\boldsymbol{e}_h|_{\boldsymbol{\beta},\mathrm{upw}}^2 \leqslant \frac{1}{2\varepsilon_3}\|\partial_t \boldsymbol{\eta}\|_{\boldsymbol{L}^2}^2 + \frac{\varepsilon_3}{2}\|\boldsymbol{e}_h\|_{\boldsymbol{L}^2}^2 \tag{55a}$$

$$+ M_{\mathrm{cnv}}\left(\frac{1}{2\varepsilon_1} + \frac{1}{\varepsilon_2}\right) |\!|\!| \boldsymbol{\eta} |\!|\!|_{\boldsymbol{\beta},\sharp}^2 + \frac{\varepsilon_1}{2}\|\boldsymbol{e}_h\|_{\boldsymbol{L}^2}^2 + \varepsilon_2 |\boldsymbol{e}_h|_{\boldsymbol{\beta},\mathrm{upw}}^2. \tag{55b}$$

Choosing $\varepsilon_1 = \varepsilon_3 = \frac{\alpha}{2} > 0$, $\varepsilon_2 = \frac{1}{2}$, and multiplication by 2 yields, for a.e. $t \in (0,T)$,

$$\frac{\mathrm{d}}{\mathrm{d}t}\|\boldsymbol{e}_h(t)\|_{\boldsymbol{L}^2}^2 + \nu C_\sigma |\!|\!| \boldsymbol{e}_h(t) |\!|\!|_e^2 + |\boldsymbol{e}_h(t)|_{\boldsymbol{\beta},\mathrm{upw}}^2 \tag{56a}$$

$$\leqslant 2\alpha^{-1}\|\partial_t \boldsymbol{\eta}\|_{\boldsymbol{L}^2}^2 + \alpha \|\boldsymbol{e}_h\|_{\boldsymbol{L}^2}^2 + M_{\mathrm{cnv}}\left(2\alpha^{-1} + 4\right) |\!|\!| \boldsymbol{\eta}(t) |\!|\!|_{\boldsymbol{\beta},\sharp}^2. \tag{56b}$$

Gronwall's lemma [21, Lemma 6.9] and $\boldsymbol{u}_h(0) = \boldsymbol{\pi}_s \boldsymbol{u}_0$ conclude the proof. ■

**Remark 3.7:** In Theorem 3.6, one can choose $\alpha = T^{-1}$ and thereby transform the exponential factor $C\alpha^{-1}e^{\alpha T}$ in the end-of-simulation time $T$ into the linear one $CT$. Note that such a weaker linear dependence is generally not possible to obtain in the nonlinear Navier–Stokes case; cf. Subsection 4.1. Alternatively to our approach, it is also possible to use a change of variables of $(\boldsymbol{u}, p)$ on the continuous level to obtain a transformed Oseen problem with an additional positive zeroth order reaction term; cf. [24]. ▲

**Remark 3.8:** Instead of using the discrete stationary Stokes projection $\boldsymbol{\pi}_s$ in the error splitting (52), one could also use the discrete Helmholtz projection; cf. [35,3]. As a consequence, in Theorem 3.6, the term involving the time derivative would vanish. However, the disadvantage is that the approximation properties of the Helmholtz projection in $\boldsymbol{H}(\mathrm{div})$ would have to be quantified. ▲

**Corollary 3.9 (*Velocity convergence rate*)**

*Under the assumptions of the previous theorem, assume a smooth solution according to*

$$\boldsymbol{u} \in L^\infty(0,T;\boldsymbol{H}^r), \quad \partial_t \boldsymbol{u} \in L^2(0,T;\boldsymbol{H}^r), \quad r > \frac{3}{2}. \tag{57}$$

*Then, with $r_{\boldsymbol{u}} = \min\{r, k+1\}$ and a constant $C$ independent of $h$ and $\nu^{-1}$, we obtain the following convergence rate:*

$$\|\boldsymbol{u} - \boldsymbol{u}_h\|_{L^\infty(\boldsymbol{L}^2)}^2 + \int_0^T \left[\nu C_\sigma |\!|\!| (\boldsymbol{u} - \boldsymbol{u}_h)(\tau) |\!|\!|_e^2 + |(\boldsymbol{u} - \boldsymbol{u}_h)(\tau)|_{\boldsymbol{\beta},\mathrm{upw}}^2\right] \mathrm{d}\tau \tag{58a}$$

$$\leqslant CTh^{2(r_{\boldsymbol{u}}-1)} \times \left[h^2\left(\|\partial_t \boldsymbol{u}\|_{L^2(\boldsymbol{H}^{r_{\boldsymbol{u}}})}^2 + \|\boldsymbol{u}\|_{L^\infty(\boldsymbol{H}^{r_{\boldsymbol{u}}})}^2\right)\right. \tag{58b}$$

$$\left. + h\|\boldsymbol{\beta}\|_\infty (\gamma + 1)\|\boldsymbol{u}\|_{L^2(\boldsymbol{H}^{r_{\boldsymbol{u}}})}^2 + \|\boldsymbol{\beta}\|_\infty \|\boldsymbol{u}\|_{L^2(\boldsymbol{H}^{r_{\boldsymbol{u}}})}^2\right] \tag{58c}$$

**Proof:** The aim is to estimate the terms on the right-hand side of Theorem 3.6. Thus, the approximation properties of the stationary Stokes projection have to be assessed. If $\boldsymbol{w}=\boldsymbol{\eta}$ and $\boldsymbol{w}=\partial_t\boldsymbol{\eta}$ is inserted in Theorem 3.3, respectively, we can directly bound $\boldsymbol{L}^2$-errors. Indeed, using that $\partial_t$ and $\boldsymbol{\pi}_s$ commute, the optimal approximation properties together with the commuting diagram property (Assumptions 3 and 4) yield

$$\|\partial_t\boldsymbol{\eta}(\tau)\|_{\boldsymbol{L}^2}^2 \leqslant C_s h^2 \inf_{\boldsymbol{v}_h\in\boldsymbol{V}_h^{\mathrm{div}}} |||\partial_t\boldsymbol{u}-\boldsymbol{v}_h|||_{e,\sharp}^2 \leqslant Ch^{2r_{\boldsymbol{u}}}|\partial_t\boldsymbol{u}|_{\boldsymbol{H}^{r_{\boldsymbol{u}}}}^2. \tag{59}$$

Details are omitted since this is a fairly standard result; cf. [19]. Proceeding, recall the $|||\cdot|||_{\boldsymbol{\beta},\sharp}$-norm from (28):

$$|||\boldsymbol{\eta}|||_{\boldsymbol{\beta},\sharp}^2 = \|\boldsymbol{\eta}\|_{\boldsymbol{L}^2}^2 + \sum_{F\in\mathcal{F}_h^i}\oint_F \frac{\gamma}{2}|\boldsymbol{\beta}\cdot\boldsymbol{n}_F||[\![\boldsymbol{\eta}]\!]|^2\,\mathrm{d}\boldsymbol{s} \tag{60a}$$

$$+\|\boldsymbol{\beta}\|_\infty^2\sum_{K\in\mathcal{T}_h} h_K^{-2}\|\boldsymbol{\eta}\|_{\boldsymbol{L}^2(K)}^2 + \|\boldsymbol{\beta}\|_\infty\sum_{K\in\mathcal{T}_h}\|\boldsymbol{\eta}\|_{\boldsymbol{L}^2(\partial K)}^2 \tag{60b}$$

The first term can be bound analogously to the time derivative, thus $\|\boldsymbol{\eta}\|_{\boldsymbol{L}^2}^2 \leqslant Ch^{2r_{\boldsymbol{u}}}|\boldsymbol{u}|_{\boldsymbol{H}^{r_{\boldsymbol{u}}}}^2$. Due to the presence of a negative power of $h$, the third term is the one which actually reduces the overall convergence order of the method. We obtain

$$\|\boldsymbol{\beta}\|_\infty^2\sum_{K\in\mathcal{T}_h} h_K^{-2}\|\boldsymbol{\eta}\|_{\boldsymbol{L}^2(K)}^2 \leqslant C\|\boldsymbol{\beta}\|_\infty^2 h^{2(r_{\boldsymbol{u}}-1)}|\boldsymbol{u}|_{\boldsymbol{H}^{r_{\boldsymbol{u}}}}^2. \tag{61}$$

For the remaining two terms, we have to use the discrete trace inequality (18). As an example we present the estimate for the first one explicitly:

$$\sum_{F\in\mathcal{F}_h^i}\oint_F \frac{\gamma}{2}|\boldsymbol{\beta}\cdot\boldsymbol{n}_F||[\![\boldsymbol{\eta}]\!]|^2\,\mathrm{d}\boldsymbol{s} \leqslant \sum_{F\in\mathcal{F}_h^i}\oint_F \gamma|\boldsymbol{\beta}\cdot\boldsymbol{n}_F|\left[|\boldsymbol{\eta}^+|^2+|\boldsymbol{\eta}^-|^2\right]\mathrm{d}\boldsymbol{s} \tag{62a}$$

$$\leqslant \|\boldsymbol{\beta}\|_\infty\sum_{K\in\mathcal{T}_h}\gamma\|\boldsymbol{\eta}\|_{\boldsymbol{L}^2(\partial K)}^2 \leqslant \|\boldsymbol{\beta}\|_\infty\sum_{K\in\mathcal{T}_h}\gamma C_{\mathrm{tr}}^2 N_\partial h_K^{-1}\|\boldsymbol{\eta}\|_{\boldsymbol{L}^2(K)}^2 \tag{62b}$$

$$\leqslant C\gamma\|\boldsymbol{\beta}\|_\infty h^{2r_{\boldsymbol{u}}-1}|\boldsymbol{u}|_{\boldsymbol{H}^{r_{\boldsymbol{u}}}}^2 \tag{62c}$$

After applying a similar argument to the last remaining term, we obtain

$$\|\boldsymbol{\beta}\|_\infty\sum_{K\in\mathcal{T}_h}\|\boldsymbol{\eta}\|_{\boldsymbol{L}^2(\partial K)}^2 \leqslant C\|\boldsymbol{\beta}\|_\infty h^{2r_{\boldsymbol{u}}-1}|\boldsymbol{u}|_{\boldsymbol{H}^{r_{\boldsymbol{u}}}}^2. \tag{63}$$

Thus,

$$\|\boldsymbol{e}_h\|_{L^\infty(\boldsymbol{L}^2)}^2 + \int_0^T\left[\nu C_\sigma|||\boldsymbol{e}_h(\tau)|||_e^2 + |\boldsymbol{e}_h(\tau)|_{\boldsymbol{\beta},\mathrm{upw}}^2\right]\mathrm{d}\tau \tag{64a}$$

$$\leqslant C\alpha^{-1}e^{\alpha T}h^{2(r_{\boldsymbol{u}}-1)}\times\Big[h^2\Big(\|\partial_t\boldsymbol{u}\|_{L^2(\boldsymbol{H}^{r_{\boldsymbol{u}}})}^2+\|\boldsymbol{u}\|_{L^2(\boldsymbol{H}^{r_{\boldsymbol{u}}})}^2\Big) \tag{64b}$$

$$+h\|\boldsymbol{\beta}\|_\infty(\gamma+1)\|\boldsymbol{u}\|_{L^2(\boldsymbol{H}^{r_{\boldsymbol{u}}})}^2 + \|\boldsymbol{\beta}\|_\infty\|\boldsymbol{u}\|_{L^2(\boldsymbol{H}^{r_{\boldsymbol{u}}})}^2\Big]. \tag{64c}$$

To finish the proof, choose $\alpha=T^{-1}$ (see Remark 3.7) and use the triangle inequality to extend this estimate to the full error. At this point, the regularity assumption of $\boldsymbol{u}$ being in $L^\infty$ in time is needed. ■

**Remark 3.10:** Provided $\boldsymbol{u}$ is sufficiently smooth, Corollary 3.9 implies that the error in the kinetic energy, the error in the energy dissipation rate and scaled jumps of the discrete velocity $\boldsymbol{u}_h$ all converge with order $k$ to zero. Thus, our $\boldsymbol{H}(\mathrm{div})$-FEM has the same asymptotical behaviour as $h \to 0$ as the analysis in [17, 24] proves for other inf-sup stable FEM. For equal-order FEM, [13] shows that one can even get half an order more. Also, (58) reveals that the limiting factor, which reduces the overall convergence to order $k$, originates from the convective term. Therefore, with $\boldsymbol{\beta} \equiv \mathbf{0}$, an order $k+1$ convergence of the $\boldsymbol{L}^2$-error for the time-dependent Stokes problem can be concluded directly. ▲

## 4 High Reynolds number numerical experiments

Our numerical experiments are conducted for the fully nonlinear Navier–Stokes equations in the high Reynolds number regime. We take advantage of the FEM package COMSOL Multiphysics 5.1. The time discretisation is performed with a fully coupled adaptive BDF(2)-scheme. In order to solve nonlinear systems, a modified Newton method with out-of-date Jacobians is used. An iteration is considered converged if its relative residual is below $10^{-6}$. Linear system are solved with the direct solver PARDISO.

All subsequent examples are two-dimensional problems with vortical structures. We always employ unstructured triangular Delaunay meshes to solve them. In order to assess the performance of $\boldsymbol{H}(\mathrm{div})$-conforming FEM, we compare our results with different FE-based numerical schemes which also make use of the following spaces (corresponding vector-valued spaces are written in bold):

$$\mathbb{P}_k = \left\{ v_h \in C(\overline{\Omega}) \colon v_h\big|_K \in \mathbb{P}_k(K),\ \forall K \in \mathcal{T}_h \right\} \tag{65a}$$

$$\mathbb{P}_k^{\mathrm{dc}} = \left\{ v_h \in L^2(\Omega) \colon v_h\big|_K \in \mathbb{P}_k(K),\ \forall K \in \mathcal{T}_h \right\} \tag{65b}$$

Then, for an $\boldsymbol{H}(\mathrm{div})$-FEM we use the Raviart–Thomas element [7]

$$\mathbb{RT}_k = \left\{ \boldsymbol{v}_h \in \boldsymbol{H}(\mathrm{div};\Omega) \colon \boldsymbol{v}_h\big|_K \in \mathbf{P}_k(K) \oplus \boldsymbol{x}\mathbb{P}_k(K),\ \forall K \in \mathcal{T}_h \right\}. \tag{66}$$

The resulting inf-sup stable FE pair is thus given by $\mathbb{RT}_k/\mathbb{P}_k^{\mathrm{dc}}$, which is abbreviated as 'RT$k$'. Due to the fact that $\nabla \cdot \mathbb{RT}_k \subseteq \mathbb{P}_k^{\mathrm{dc}}$ (cf. [7, Proposition 2.3.3]), Assumption A1 is fulfilled and we obtain an exactly divergence-free method. For the jump penalisation parameter, $\sigma = 6\frac{(k+1)(k+d)}{d}$ is chosen in all examples. Here, we follow the asymptotic behaviour w.r.t. the polynomial degree $k \geqslant 1$ suggested by [28]. However, we found that the results are rather insensitive towards this parameter as long as it is large enough to guarantee discrete coercivity.

Moreover, as a divergence-free $\boldsymbol{H}^1$-conforming method, we employ the Scott–Vogelius element [43] with velocity/pressure space $\mathbf{P}_k/\mathbb{P}_{k-1}^{\mathrm{dc}}$. Inf-sup stability is guaranteed by using barycentre-refined meshes. This method is denoted 'SV$k$'; for more information, we refer to the authors' previous work [42].

As a representative of inf-sup stable $\boldsymbol{H}^1$-conforming methods which are not exactly divergence-free, we take the well-known Taylor–Hood method of order $k$ with velocity/pressure pair $\mathbf{P}_k/\mathbb{P}_{k-1}$. The nonlinear term is treated with the EMAC formulation [14]. Therein, as well as in the authors' work [42], it is shown that the EMAC formulation holds several theoretical and practical advantages over more common formulations as, for example, the convective or skew-symmetric formulation. This method is abbreviated as 'eTH$k$'. Nevertheless, we also tried the convective and skew-symmetric formulations. However, the EMAC formulation is never inferior to them (mostly, it yields significantly better results) and therefore, we restrict ourselves to showing results exclusively for the EMAC formulation.

Furthermore, we want to consider an $\boldsymbol{H}^1$-conforming equal-order method (not exactly divergence-free), and thus choose the FE pair $\mathbb{P}_k/\mathbb{P}_k$. For the convective part, as is probably most common, a skew-symmetric formulation is applied and we dub the resulting method 'EO$k$'. Note that, in order to obtain a stable method with equal-order interpolation, the pressure always has to be stabilised; cf. [12]. In accordance with [11], we choose the parameter $\gamma_{\mathrm{PS}} = 0.01$ and, for all applications of EO$k$, add the term

$$\gamma_{\mathrm{PS}} \sum_{F\in\mathcal{F}_h^i} \oint_F h_F^2 [\![\nabla p_h]\!] \cdot [\![\nabla q_h]\!] \,\mathrm{d}\boldsymbol{s}. \tag{67}$$

For all above introduced FE schemes, additionally to the corresponding 'Galerkin' formulations (those formulations which guarantee stability), we also consider the possibility of adding suitable stabilisation terms. For non-divergence-free methods, the divergence of discrete velocities is stabilised using a grad-div term (GD) with parameter $\gamma_{\mathrm{GD}} \geqslant 0$ [31]

$$\gamma_{\mathrm{GD}} \sum_{K\in\mathcal{T}_h} \int_K (\nabla\cdot\boldsymbol{u}_h)(\nabla\cdot\boldsymbol{v}_h)\,\mathrm{d}\boldsymbol{x}. \tag{68}$$

Whenever we add grad-div stabilisation the parameter $\gamma_{\mathrm{GD}} = 0.01$ is chosen. Choosing $\gamma_{\mathrm{GD}}$ as a constant is of course not necessary. For a more elaborate discussion on the choice of the grad-div parameter and associated difficulties and specialities we refer to [29,1] where stationary incompressible flow problems are considered. For divergence-free methods, grad-div stabilisation is of course superfluous. However, it might be advantageous to explicitly include some (additional) kind of convection stabilisation in the numerical schemes. Note that in the $\boldsymbol{H}(\mathrm{div})$-FEM, the natural upwind stabilisation is already included as a form of convection stabilisation. From the plethora of available stabilisations we choose the gradient jump stabilisation which is usually considered in the context of continuous interior penalty (CIP) methods with parameter $\gamma_{\mathrm{CIP}} \geqslant 0$; cf., for example, [12]:

$$\gamma_{\mathrm{CIP}} \sum_{F\in\mathcal{F}_h^i} \oint_F h_F^2 |\boldsymbol{u}_h\cdot\boldsymbol{n}_F|^2 [\![\nabla\boldsymbol{u}_h]\!] : [\![\nabla\boldsymbol{v}_h]\!] \,\mathrm{d}\boldsymbol{s} \tag{69}$$

CIP stabilisation can be added to eTH$k$, EO$k$, SV$k$ and RT$k$ since all methods yield globally discontinuous discrete velocity gradients. Whenever we add it, $\gamma_{\mathrm{CIP}} = 0.1$ is chosen.

Results are often compared in terms of kinetic energy $\mathcal{K}$ and enstrophy $\mathcal{E}$. For a 2D velocity $\boldsymbol{w} = (w_1,w_2)^\dagger$, we agree on $\boldsymbol{\omega} = \nabla\times\boldsymbol{w} = \partial_{x_1} w_2 - \partial_{x_2} w_1$ and use the following:

$$\mathcal{K}(\boldsymbol{w},t) = \frac{1}{2}\,\|\boldsymbol{w}(t)\|_{\boldsymbol{L}^2}^2 = \frac{1}{2}\int_\Omega |\boldsymbol{w}(t,\boldsymbol{x})|^2\,\mathrm{d}\boldsymbol{x} \tag{70a}$$

$$\mathcal{E}(\boldsymbol{w},t) = \frac{1}{2}\,\|\nabla_h\times\boldsymbol{w}(t)\|_{L^2}^2 = \frac{1}{2}\int_\Omega |\nabla_h\times\boldsymbol{w}(t,\boldsymbol{x})|^2\,\mathrm{d}\boldsymbol{x} \tag{70b}$$

In the remainder of this section, our approach is as follows. We concentrate on 2D problems exclusively because we believe that it is important to first understand how a numerical method performs in this situation. If it does not work for 2D, there is no real possibility (or hope) that it will work satisfactorily in 3D. Our first example is the planar lattice flow for which an exact solution is known. For this problem we compare the above introduced methods and try to single out something like the 'best' $\boldsymbol{H}(\operatorname{div})$ method for this problem. Then, we consider the Kelvin–Helmholtz instabilities triggered by a mixing layer. Here, we only consider one $\boldsymbol{H}(\operatorname{div})$-FEM but show how sensitive the solution is towards mesh refinement. Finally, we apply this method on one fixed mesh to the simulation of freely decaying 2D turbulence. We include aspects of energy and enstrophy, but also of the self-organisation into large-scale structures involving energy spectra.

## 4.1 Planar lattice flow

In this section we consider the evolution of an initial velocity, which solves the stationary incompressible Euler equation, in a viscous incompressible Navier–Stokes flow. This example has also been investigated in detail in [42] and is called 'planar lattice flow' [5]. For $\boldsymbol{x} \in \Omega = (0,1)^2$, both the initial condition $\boldsymbol{u}_0$ and the corresponding known exact velocity/pressure pair $(\boldsymbol{u}, p)$ for $\nu \geqslant 0$ are given as follows:

$$\boldsymbol{u}_0(\boldsymbol{x}) = \begin{bmatrix} \sin(2\pi x_1)\sin(2\pi x_2) \\ \cos(2\pi x_1)\cos(2\pi x_2) \end{bmatrix}, \quad \boldsymbol{u}(t,\boldsymbol{x}) = \boldsymbol{u}_0(\boldsymbol{x}) e^{-8\pi^2 \nu t} \tag{71a}$$

$$p(t,\boldsymbol{x}) = \frac{1}{4}\left[\cos(4\pi x_1) - \cos(4\pi x_2)\right] e^{-16\pi^2 \nu t} \tag{71b}$$

Here, the initial velocity $\boldsymbol{u}_0$ induces a flow structure which, due to its saddle point character, is 'dynamically unstable so that small perturbations result in a very chaotic motion' [36]. We impose periodic boundary conditions on the vertical and horizontal walls of $\partial\Omega$, respectively, and the integral zero-mean condition is imposed on the pressure. There is no external forcing in this problem, that is, $\boldsymbol{f} = \boldsymbol{0}$, and $\nu = 4 \times 10^{-6}$ is fixed. For a more qualitative approach to this problem, we refer to [42], where one can get a better feeling for the appearance and behaviour of this particular flow problem.

In Table 1, an overview of the meshes and DOFs for this problem is given. Especially, the number of DOFs is split based on how many are used for the velocity and pressure discretisation, respectively. Note that for SV$k$, unstructured barycentre-refined Delaunay meshes are used while eTH$k$, EO$k$ and RT$k$ are based on unstructured Delaunay triangulations. The use of unstructured meshes introduces an additional difficulty because it makes it harder for any numerical method to preserve the symmetric nature of the flow. The meshes are chosen in such a way that the total number of DOFs for all methods approximately coincides. It is interesting to acknowledge the different distribution of DOFs for velocity and pressure for the different schemes. Especially, the SV$k$ methods spend a lot for the pressure, which can be considered as a disadvantage in problems where the pressure is not of primary interest.

Let us begin with the comparison of results for second-order FE pairs. In Figure 1, the evolution of errors w.r.t. kinetic energy $\mathcal{K}$ and enstrophy $\mathcal{E}$ can be seen for the particular 'Galerkin' formulations and for some suitably stabilised variants thereof. The first apparent conclusion is that, for each method, the stabilised variant significantly outperforms its basic, stable counterpart. For the 'Galerkin' methods, eTH2 yields worse results than SV2 which, in turn, is inferior to EO2, which is not as good as RT2. After adding suitable stabilisation, however, eTH2 and SV2 now yield comparable results. The potential advantage of a div-free $\boldsymbol{H}^1$-conforming method is thus relativised if the EMAC method is equipped with sufficient additional stabilisation. Furthermore, we observe that after adding CIP and grad-div stabilisation to the $\boldsymbol{H}^1$-conforming EO2 method, RT2 and EO2 now yield comparable results. Here, the $\boldsymbol{H}(\mathrm{div})$-conforming RT2 method benefits from the possibility to include upwinding. At this point, we want to mention that equal-order methods EO$k$ turned out to be relatively sensible towards the stabilisation parameters (especially the pressure gradient jump penalisation parameter $\gamma_{\mathrm{PS}}$), which is not very attractive from the application-oriented perspective. The main conclusion from Figure 1 is that we can discard the pure 'Galerkin' methods and instead exclusively concentrate on suitably stabilised schemes.

**Table 1** Overview of meshes and DOFs for all methods which are compared for the planar lattice flow. Abbreviations of different methods: Non-div-free $\boldsymbol{H}^1$ EMAC Taylor–Hood (eTH$k$) and equal-order (EO$k$), div-free $\boldsymbol{H}^1$ Scott–Vogelius (SV$k$) and div-free $\boldsymbol{H}(\mathrm{div})$ Raviart–Thomas (RT$k$).

| Method | $\max_{K\in\mathcal{T}_h} h_K$ | $\min_{K\in\mathcal{T}_h} h_K$ | #{triangles} | #{$\boldsymbol{u}$ DOFs} | #{$p$ DOFs} | #{DOFs} |
|---|---|---|---|---|---|---|
| RT2 | 0.0267 | 0.0104 | 8414 | 88 617 | 50 484 | 139 101 |
| SV2 | 0.0292 | 0.0087 | 20 574 | 82 618 | 61 722 | 144 340 |
| eTH2 | 0.0134 | 0.0048 | 35 204 | 141 546 | 17 785 | 159 331 |
| EO2 | 0.0158 | 0.0058 | 25 400 | 102 218 | 51 109 | 153 327 |
| RT3 | 0.0583 | 0.0247 | 1680 | 30 400 | 16 800 | 47 200 |
| SV3 | 0.0667 | 0.0221 | 3930 | 35 588 | 23 580 | 59 168 |
| eTH3 | 0.0319 | 0.0140 | 5586 | 50 720 | 11 321 | 62 041 |
| EO3 | 0.0368 | 0.0136 | 4080 | 37 094 | 18 547 | 55 641 |

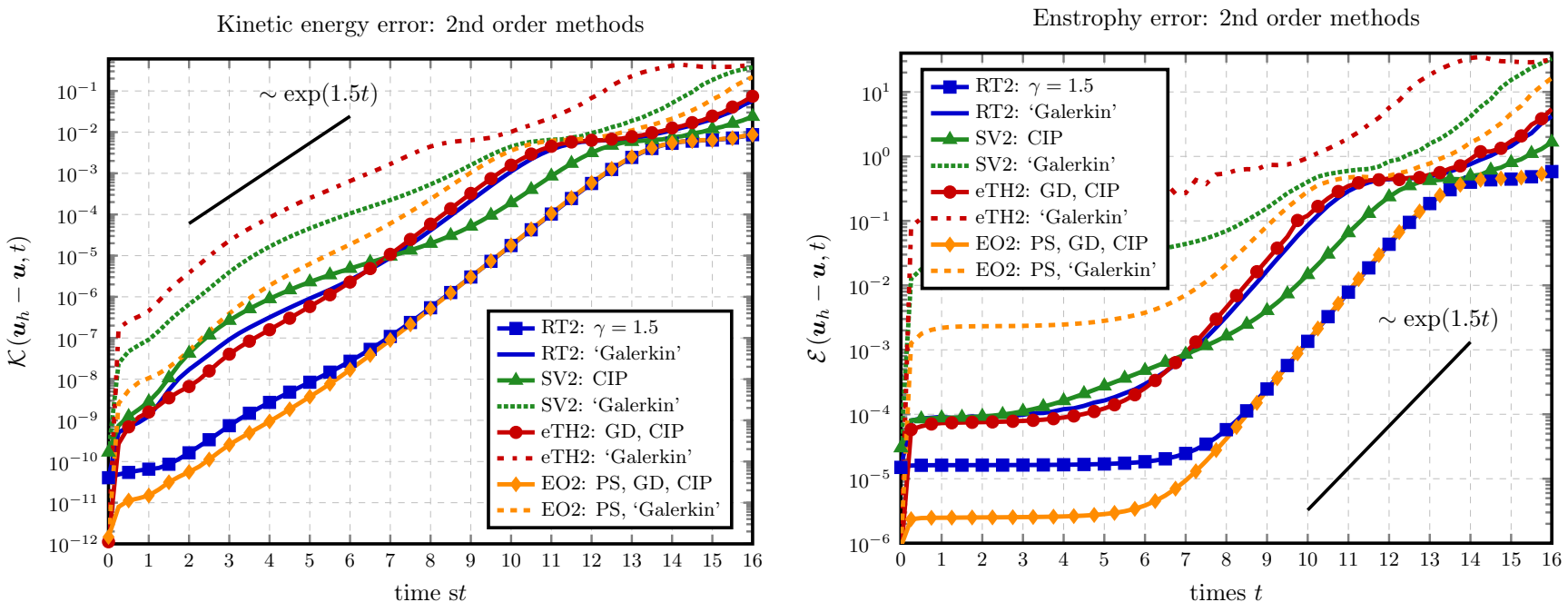

**Fig. 1** Evolution of errors w.r.t. kinetic energy $\mathcal{K}$ and enstrophy $\mathcal{E}$ for 2nd order methods with $\nu = 4 \times 10^{-6}$. The mesh data and DOFs are displayed in Table 1. $\gamma$ is the upwind parameter and PS denotes pressure stabilisation. For the other parameters, $\gamma_{\mathrm{CIP}} = 0.1$ and $\gamma_{\mathrm{GD}} = 0.01$ are chosen. 'Galerkin' denotes the stable basic variant of each method, respectively.

In Figure 2, results for third-order FE pairs are shown. In view of Table 1, it should be noted that even though the 3rd order methods use considerably less DOFs compared to the 2nd order methods, for the inf-sup stable methods, the total errors are still significantly smaller. This is a strong argument for using higher-order methods. Only for the equal-order methods, going from second to third-order, while coarsening the mesh, does not yield a considerable improvement. Returning to Figure 2, we observe that the overall tendency observed for the second-order methods can be updated: RT3 yields the best results, followed closely by EO3, and SV3 and eTH3 are roughly on par with each other. Next, the possibility of changing the convection stabilisation for RT3 is explored. To begin with, we add the same CIP gradient jump penalisation term to RT3 which improves all $\boldsymbol{H}^1$-conforming methods. The solid blue line indicates that this kind of stabilisation corrupts the $\boldsymbol{H}(\mathrm{div})$-conforming method significantly. Another possibility is to adjust the upwind parameter $\gamma$. However, the dotted and dashed blue lines show that neither increasing nor decreasing $\gamma$ away from the standard value $\gamma = 1.5$ seems to be promising. Thus, RT3 with $\gamma = 1.5$ yields very convincing results for the planar lattice flow problem which are better than the ones obtained by comparable schemes. In addition, the effort of fine-tuning stabilisation parameters can be minimised with such a method—upwinding represents a natural stabilisation which, very importantly, is not gradient-based. Moreover, numerical experiments revealed that in practice, an upwind parameter between 1 and 2 always leads to good results.

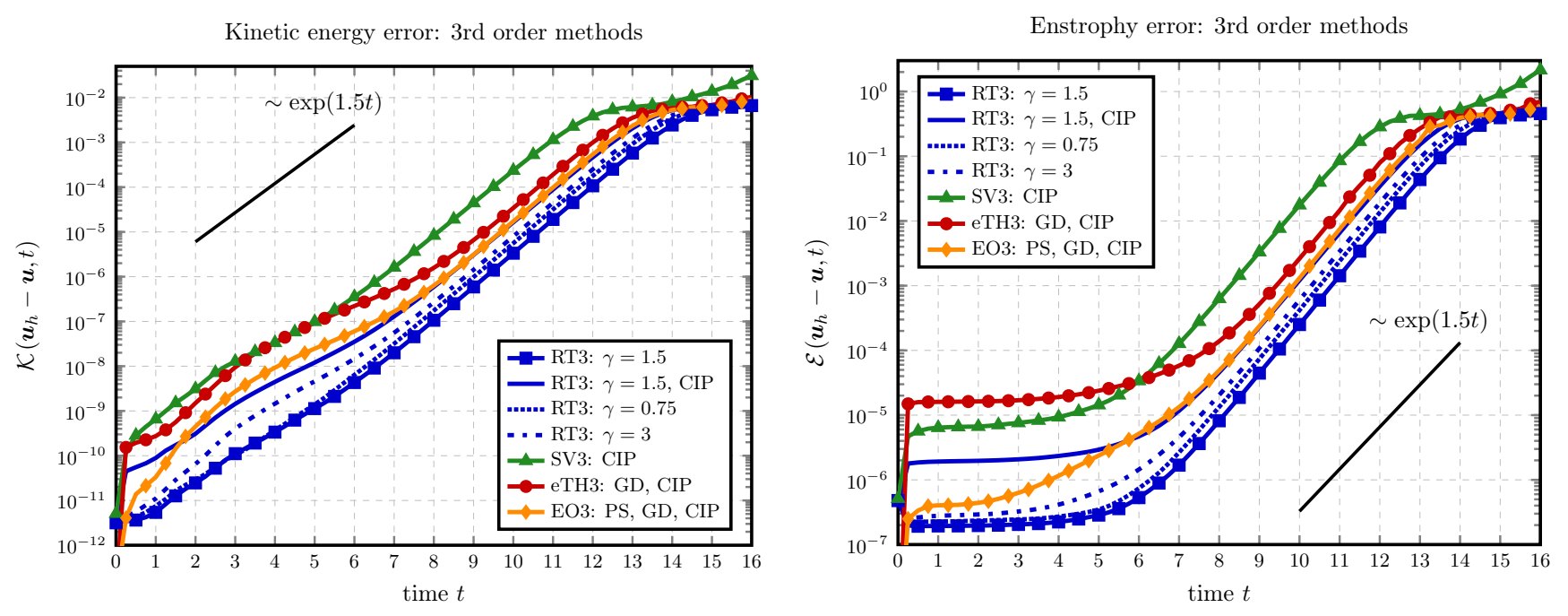

**Fig. 2** Evolution of errors w.r.t. kinetic energy $\mathcal{K}$ and enstrophy $\mathcal{E}$ for 3rd order methods with $\nu = 4 \times 10^{-6}$. The mesh data and DOFs are displayed in Table 1. $\gamma$ is the upwind parameter, PS denotes pressure stabilisation and $\gamma_{\mathrm{CIP}} = 0.1$, $\gamma_{\mathrm{GD}} = 0.01$.

Lastly, we want to take a closer look at the quality of the pressure approximation for each particular method. As can be seen in Table 1, the amount of DOFs spend for the pressure varies widely. The Scott–Vogelius methods have the most pressure DOFs, followed by the equal-order and Raviart–Thomas methods. A discretisation based on Taylor–Hood elements is rather inexpensive in terms of pressure DOFs. However, whilst eTH$k$ and SV$k$ use piecewise polynomials of order $k-1$, the corresponding pressure spaces for RT$k$ and EO$k$ consist of piecewise $k$-th order polynomials. In Figure 3, the evolution of the pressure errors for all considered methods is shown. We conclude that for 2nd order methods, EO2 and RT2 are again comparable. However, RT3 ($\gamma = 1.5$) yields the best pressure approximation for the third-order methods. Thus, even though SV$k$ always uses much more DOFs, this additional cost does not seem to pay off in terms of accuracy.

*Intermediate summary:* The application of different methods for the planar lattice flow revealed the following insights. Independent of which method is chosen, a suitable stabilisation always improves the method compared with the particular stable basic ('Galerkin') variant. Gradient-based CIP convection stabilisation corrupts $\boldsymbol{H}(\mathrm{div})$-conforming methods, whereas all $\boldsymbol{H}^1$-conforming discretisations benefit from it. Going to higher-order approximations is always worthwhile because, all in all, more accurate results can be achieved with less DOFs. For the equal-order methods, this improvement is not as pronounced as for the inherently inf-sup stable ones. Third-order div-free $\boldsymbol{H}(\mathrm{div})$-FEMs with upwinding yield the best results. Stabilised equal-order FEM ($\boldsymbol{H}^1$, non-div-free) are generally only slightly inferior. Considerably inferior are stabilised Scott–Vogelius ($\boldsymbol{H}^1$, div-free) and stabilised EMAC Taylor–Hood ($\boldsymbol{H}^1$, non-div-free) methods, which, in turn, are roughly of equal value. Remarkably, all above shown error plots clearly show an $\exp(Ct)$ behaviour (see the black asymptotes) for the Navier–Stokes problem. This is in contrast to the Oseen problem, where it is possible to prove a linear in $t$ growth of the error; cf. Remark 3.7.

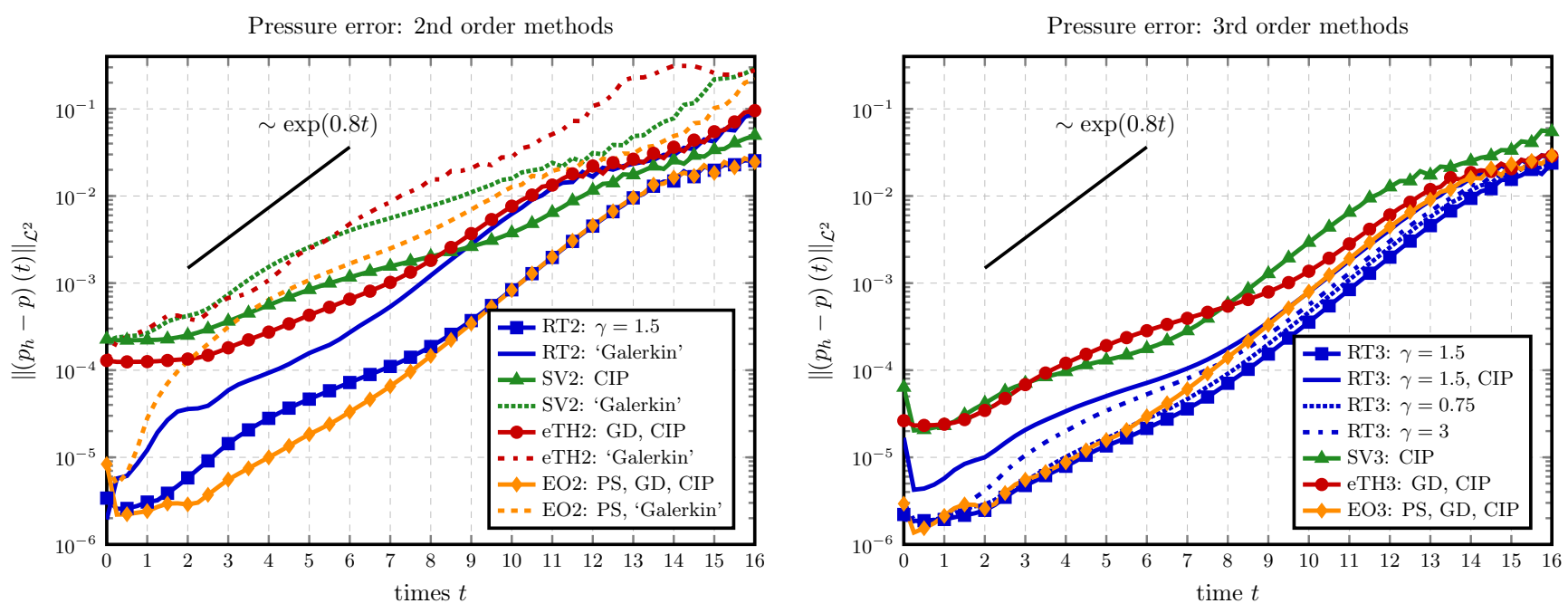

**Fig. 3** Evolution of pressure errors for second- and third-order methods with $\nu = 4 \times 10^{-6}$.

### 4.2 2D Kelvin–Helmholtz instability

From the point of view of practical applications, we now turn to a more relevant example. Two-dimensional Kelvin–Helmholtz instabilities of plane mixing layers are important test cases in fluid dynamics. Even though this kind of flow does not lead to what is typically called 'turbulence', it is extremely sensitive towards initial conditions and shows an energy spectrum $E(\kappa)$ of slope between $\kappa^{-3}$ and $\kappa^{-4}$ [34].

Based on [25,11,2], we briefly summarise the setting of the problem. On $\Omega = (0,1)^2$, periodic boundary conditions for $x_1 \in \{0,1\}$ and free-slip boundary conditions at $x_2 \in \{0,1\}$ are imposed. There is no external forcing, that is, $\boldsymbol{f} = \boldsymbol{0}$. Similarly to the planar lattice flow, the whole problem is determined by an initial condition which evolves in a viscous incompressible flow. However, the behaviour of the flow will be more dynamic. Let $\delta_0 = 1/28$ denote the initial vorticity thickness, $u_\infty = 1$ be a reference velocity, $c_n = 10^{-3}$ define a scaling factor and choose the viscosity according to $\nu^{-1} = 28 \times 10^4$. Thus, the Reynolds number associated with this problem is $Re = u_\infty \delta_0 \nu^{-1} = 10^4$. Introducing the stream function

$$\psi(x_1, x_2) = c_n u_\infty \exp\left(-\frac{(x_2 - 0.5)^2}{\delta_0^2}\right)\left[\cos(8\pi x_1) + \cos(20\pi x_1)\right], \tag{72}$$

the initial velocity field for our simulation is given as follows:

$$\boldsymbol{u}_0(\boldsymbol{x}) = \begin{bmatrix} u_\infty \tanh\left(\frac{2x_2 - 1}{\delta_0}\right) \\ 0 \end{bmatrix} + \begin{bmatrix} \partial_{x_2} \psi \\ -\partial_{x_1} \psi \end{bmatrix} \tag{73}$$

**Table 2** Overview of meshes and DOFs for the 2D Kelvin–Helmholtz instability.

| Name | $\max_{K\in\mathcal{T}_h} h_K$ | #{triangles} | #{$\boldsymbol{u}$ DOFs} | #{$p$ DOFs} |
|---|---|---|---|---|
| RT3-a | $2.425\times10^{-2}$ | 10 602 | 191 236 | 106 020 |
| RT3-b | $1.159\times10^{-2}$ | 47 646 | 858 476 | 476 460 |
| RT3-c | $7.487\times10^{-3}$ | 119 602 | 2 154 172 | 1 196 020 |
| RT3-d | $4.924\times10^{-3}$ | 268 762 | 4 839 716 | 2 687 620 |

For the evaluation, the scaled time unit $\bar{t} = \delta_0/u_\infty$ is introduced. All simulations are computed up until $T = 7.1429 = 200\bar{t}$. Moreover, in the context of this mixing layer problem, it is an established procedure to consider the following vorticity thickness $\delta(t)$:

$$\delta(t) = \frac{2u_\infty}{\omega_{\max}(t)}, \qquad \omega_{\max}(t) = \sup_{x_2\in[0,1]} |\langle\omega\rangle(t,x_2)| = \sup_{x_2\in[0,1]} \left| \int_0^1 \omega(t,\boldsymbol{x})\,\mathrm{d}x_1 \right| \tag{74}$$

In practice, however, the spatial supremum in $\omega_{\max}$ is always approximated by taking the maximum absolute value of the integral over 64 equidistantly spaced lines parallel to the $x_1$-axis. Note that these evaluation lines do not coincide with mesh lines since unstructured Delaunay triangulations are used. Moreover, in our experience, using more than 64 lines does not noticeably change the resulting vorticity thickness. In time, $\delta(t)$ is evaluated 200 times for $t \in \{1,\dots,200\}\times\bar{t}$ which yields a sufficiently resolved evolution.

In Section 4.1 we compared various different FEM for a problem with an analytical solution. As a conclusion, we singled out that the third-order $\boldsymbol{H}(\mathrm{div})$-conforming and divergence-free FEM RT3 with upwind stabilisation ($\gamma = 1.5$) gives very convincing results. Thus, we exclusively use this method for the simulation of the 2D Kelvin–Helmholtz instability. However, we compare different levels of resolution by employing RT3 on a sequence of meshes which represent under-resolved to reasonably well-resolved situations; see Table 2. Note that all methods in [25, 11, 2] are 2nd order. Moreover, the particular number of velocity DOFs in those references is always about 100 000 and thus comparable at most to our coarsest mesh. Better resolving simulations can be found, for example, in [30] where a second-order method is used also on a mesh with $h \approx 5\times10^{-3}$ (note that the problem considered in [30] is not comparable quantitatively since it considers a larger initial vorticity thickness of $1/14$). However, to the best of our knowledge, with our third-order simulations this work presents the most resolved results for this 2D Kelvin–Helmholtz instability in the literature.

Prior to a more quantitative analysis, let us first understand the general behaviour of the flow and compare our results to [25, 11, 2]. To this end, the evolution of the vorticity $\nabla\times\boldsymbol{u}_h$, obtained with RT3-d, can be seen in Figure 4. We draw the following conclusions:

- *Four primary vortices and their pairing:* In agreement with the other references, 4 primary vortices develop between 10 and 20 time units $\bar{t}$. These vortices merge at about $35\bar{t}$ which is also observed in [25, 2]. In [11] this pairing takes place later.
- *Pairing of two secondary vortices:* The two secondary vortices are standing for a certain amount of time. However, the instance in time where the second pairing begins is strongly dependent on which method and resolution is used. For example, considering $t = 100\bar{t}$, our two primary vortices in Figure 4 are still clearly separated and are aligned on a line parallel to the $x_1$-axis. In contrast, the two vortices in [25, 2] are already moving towards each other near the periodic boundary. In [11], the pairing has already begun but the last vortex seems to be located in the center of the domain.

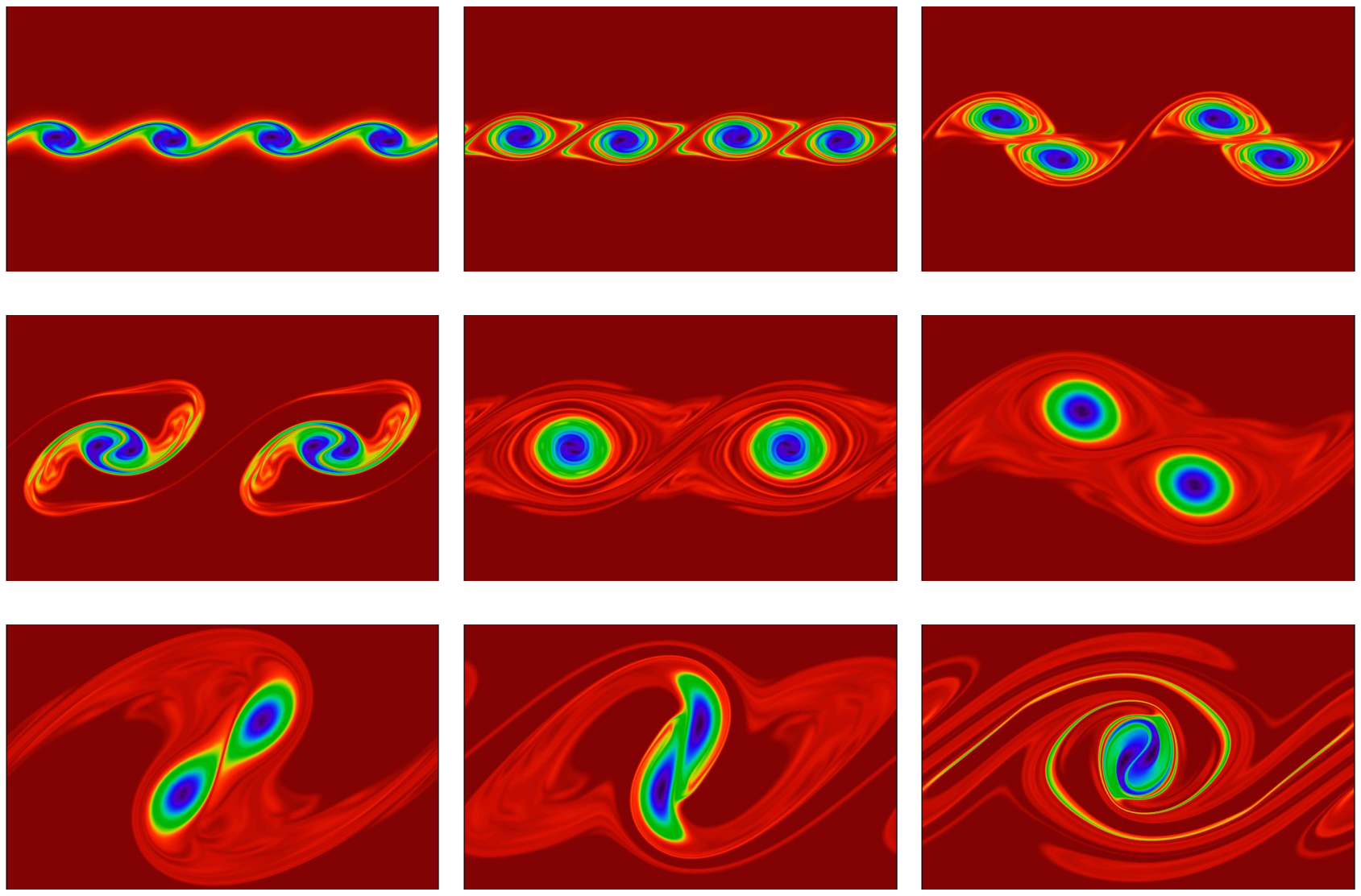

**Fig. 4** Vorticity $\nabla \times \boldsymbol{u}_h(t)$ for 2D Kelvin–Helmholtz instability at (from left to right and top to bottom) $t \in \{10, 20, 30, 40, 100, 155, 165, 180, 200\} \times \bar{t}$. Obtained with div-free $\boldsymbol{H}(\mathrm{div})$-FEM RT3-d; cf. Table 2.

– *Position of last rotating vortex:* Independent of when the last pairing occurs, there is no consensus concerning the location of the last vortex. Our results here, as well as the ones presented in [11], support the claim that the last vortex should rotate in the centre of the domain. In [25,2], the last rotation takes place across the periodic boundary.

Moreover, we want to draw attention to the fact that between the main vortices, fine-scale flow structures can be observed very well. Such structures are not dissipated numerically by the $\boldsymbol{H}(\mathrm{div})$-FEM. For a more detailed description of the mechanisms behind vortex merging, we refer to [37].

Furthermore, the evolution of both kinetic energy $\mathcal{K}$ and enstrophy $\mathcal{E}$ can be seen in Figure 5. Roughly independent of the mesh size, the kinetic energy decays only very slowly in our simulations (decay in energy is about 0.3 % ). This is in agreement with [11] but in contrast to [25] (about 1 % energy loss) and [2] (about 5 % energy loss). We interpret this observation as an indicator that $\boldsymbol{H}(\mathrm{div})$-FEM do have a much less dissipative nature (even on coarse meshes) compared to other methods. We conjecture that the main reason for this behaviour lies in the minimum amount of stabilisation (only upwinding) which is needed for $\boldsymbol{H}(\mathrm{div})$-conforming FEM. Regarding the evolution of enstrophy we observe that a more accurate method with a higher resolution leads to a later decrease in enstrophy. Actually, the stages of the enstrophy are directly connected to the pairing of vortices in the simulation. Especially the occurrence of the last pairing is very clearly observable by the sudden decrease in enstrophy towards the end of the simulations.

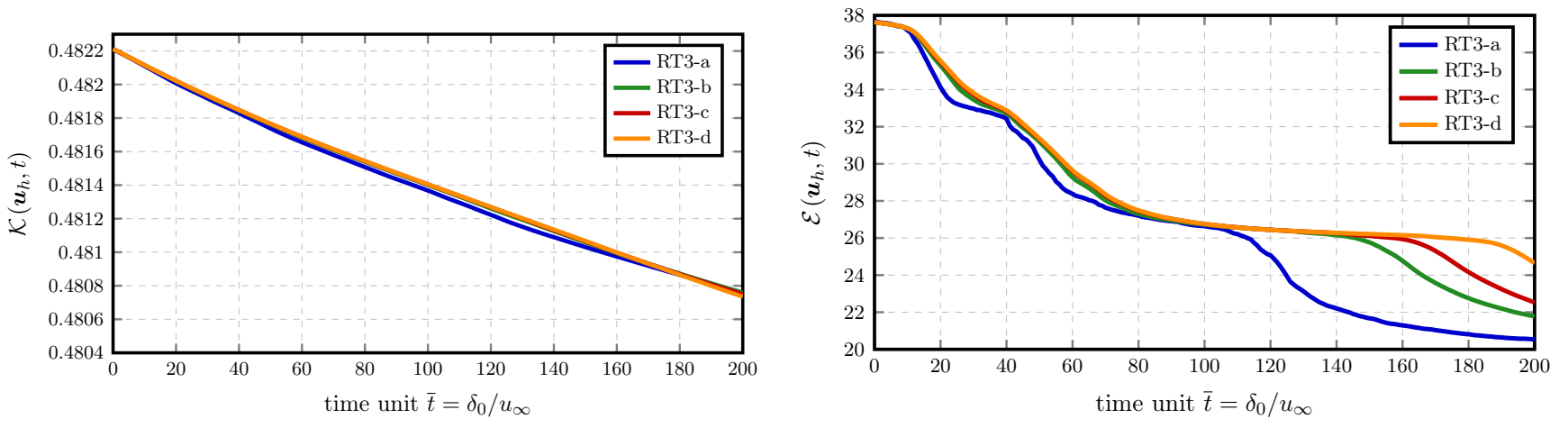

**Fig. 5** Evolution of kinetic energy $\mathcal{K}$ (left) and enstrophy $\mathcal{E}$ (right); cf. (70).

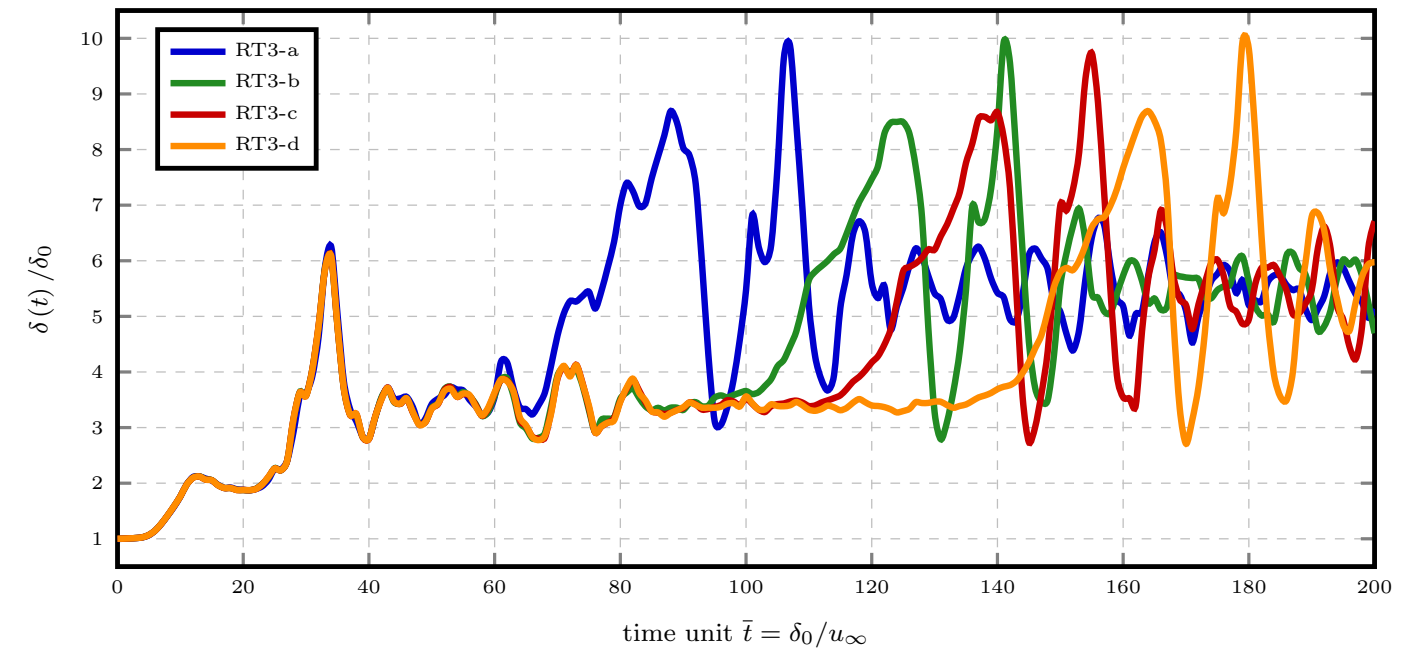

**Fig. 6** Evolution of scaled vorticity thickness $\delta(t)/\delta_0$; cf. (74).

In Figure 6, the scaled vorticity thickness is shown. Mainly, this plot illustrates in more detail when the pairing processes occur in time. The general tendency that a well-resolved simulation tends to preserve multiple vortices as long as possible is again clearly reflected. Particularly the last pairing process, where two vortices merge to become one, is very sensible with respect to how accurate the simulation is. Despite a temporal shift of the last pairing compared to [25], however, the total values and the general behaviour of our scaled vorticity thickness agrees quite well. Whilst our simulations indicate a maximum over time of $\delta(t)/\delta_0$ of about 10, the results in [2] show a higher maximum value of about 12.

A likely reason for the mentioned time lag of the last pairing is the above discussed difference w.r.t. diffusivity (or dissipativity) of a numerical method. It is much easier for a more diffusive numerical method to show behaviour resembling mesh convergence since the fine scales of the flow, which would normally need a higher resolution to be resolved properly, are smoothed out by means of numerical diffusion. Therefore, our simulations of the 2D Kelvin–Helmholtz instability problem can be understood as another example which emphasises the importance of approximation schemes with as little as possible numerical smearing.

## 4.3 Freely decaying 2D turbulence

In this section we combine both the planar lattice flow problem and the 2D Kelvin–Helmholtz instability and extend them to obtain an even more difficult situation. The original idea in Section 4.1 was to place four oppositely rotating vortices into a periodic square domain at $t = 0$ and study their evolution in absence of external forcing and boundary conditions. In doing so, the main question was to assess the ability of numerical schemes to preserve the structure of the flow, given by means of the initial condition—a good numerical method gives four standing vortices which decay in a stable manner over time. In Section 4.2, on the other hand, an inherently dynamic problem arises due to an initially imposed shear layer and merging of co-rotating vortices (like-signed vorticity regions) can be observed.

In this new example, however, the initial velocity $\boldsymbol{u}_0 = (\partial_{x_2}\psi, -\partial_{x_1}\psi)^\dagger$ shall represent $n_v^2 = 32^2 = 1024$ pairwise oppositely rotating vortices resulting from the stream function

$$\psi(\boldsymbol{x}) = 10^{-2}\sum_{k,j=1}^{n_v}(-1)^{k+j}\exp\left(-10^4\left[\left(x_1 - \frac{k}{n_v+1}\right)^2 + \left(x_2 - \frac{j}{n_v+1}\right)^2\right]\right). \tag{75}$$

Here, $\boldsymbol{u}_0$ evolves unimpeded; thus $\boldsymbol{f} = \boldsymbol{0}$ and therefore we have a freely decaying problem. The domain is $\Omega = (0,1)^2$ and periodic boundary conditions are imposed on the vertical and horizontal walls of $\partial\Omega$, respectively. Thus, the integral zero-mean condition is required for the pressure. Fully intentionally, we do not distribute the vortices equidistantly—the distance of vortices across the periodic boundaries is greater than in the 'interior' of the domain.

A flow with such an initial condition, especially for high Reynolds numbers, is very unstable and tends to evolve into a rather chaotic motion. This phenomenon is known as two-dimensional turbulence; we refer to [44, 6] for more information. 2D turbulence follows the Kraichnan–Batchelor–Leith (KBL) theory [18] and typical properties of freely decaying flows are the energy spectrum $E(\kappa) \sim \kappa^{-3}$ and, in stark contrast to 3D turbulence, the self-organisation of small-scale features of the flow into constantly growing large-scale coherent vortices; cf. [26]. The last aspect is connected to the presence of an additional inverse cascade in 2D turbulence. Studying such problems is not novel; see, for example, [10] where freely decaying 2D turbulence has been studied in the vorticity-stream function formulation. Also, we would like to mention [40] where a comparison of different numerical schemes (no FEM, though) for the DNS of freely decaying 2D turbulence in $\omega/\psi$-formulation has been presented. Also, in the context of atmospheric flows, the transfer of energy and enstrophy between scales is very important and a comparable flow configuration for the study of freely decaying 2D turbulence can be found in [45] (also in $\omega/\psi$-formulation).

However, to the best of the authors' knowledge, in the literature, there are no comparable studies available for the original velocity-pressure formulation. Thus, we want to fill this gap by showing that our $\boldsymbol{H}(\mathrm{div})$-FEM is able to produce trustworthy simulations for freely decaying 2D turbulence. To this end, let us consider simulations with the three different viscosities $\nu \in \left\{5\times10^{-5}, 10^{-5}, 4\times10^{-6}\right\}$. Concerning the mesh resolution, the Kelvin–Helmholtz instability in Section 4.2 already revealed that it is extremely expensive to obtain mesh-converged solutions. Therefore, for this example, we restrict ourselves to exclusively one unstructured Delaunay mesh consisting of 119 602 triangles. With our favourite RT3 method ($\gamma = 1.5$), this leads to 2 154 172 velocity and 1 196 020 pressure DOFs.

First of all, Figure 7 shows the evolution of the kinetic energy and enstrophy for the three different viscosities. We observe the expected behaviour that with smaller $\nu$, $\mathcal{K}$ decays much more slowly. More interesting is the behaviour of the enstrophy. Especially for the two smaller viscosities, one can observe a small initial range up to $t \approx 0.75$ where the enstrophy decays only slowly. In accordance with [18, Chapter 10] this describes a transition zone in which 2D turbulence develops out of the ordered initial condition. After that, a stronger decay in $\mathcal{E}$ can be observed which corresponds to fully developed 2D turbulence.

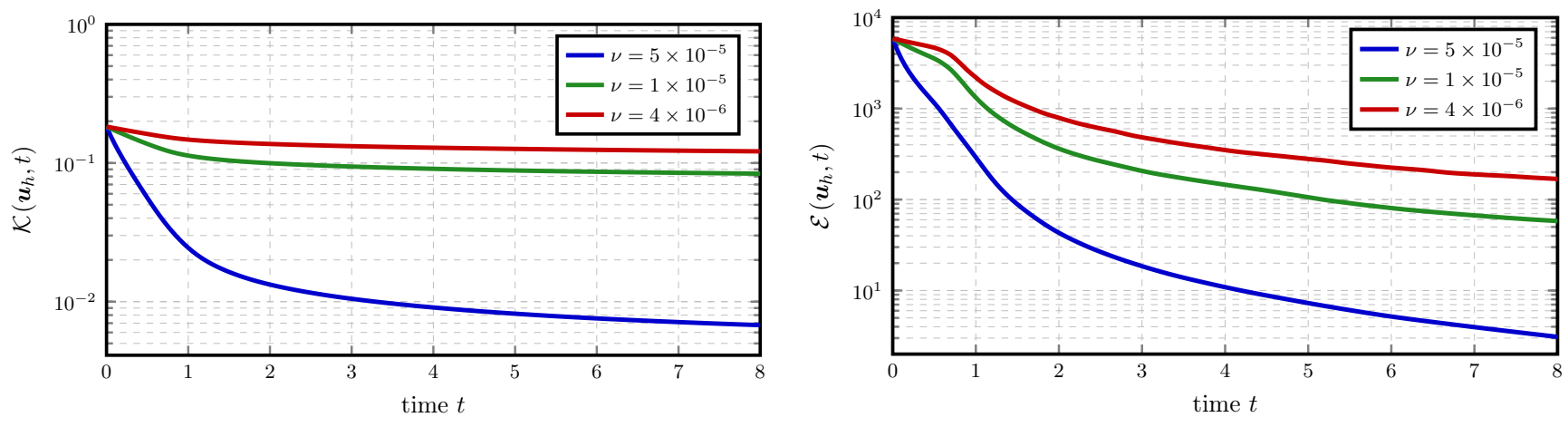

**Fig. 7** Evolution of kinetic energy $\mathcal{K}$ (left) and enstrophy $\mathcal{E}$ (right) for $\nu \in \{5 \times 10^{-5}, 10^{-5}, 4 \times 10^{-6}\}$.

In Figure 8, snapshots of the vorticity can be seen for the different viscosities (columns) at certain time instances. At first, it becomes clear how viscous forces attack the initial vorticity field, consisting of 1024 clearly separated vortices, and therefore, depending on which $\nu$ is considered, a more or less chaotic motion can be observed. Here, a smaller viscosity leads to a more small-scale structure of the flow. As time proceeds, one can directly observe that while the multitude of vortices moves trough the domain, like-signed vortices merge and oppositely rotating vortices repel each other. One can see clearly that over time, more and more large-scale structures develop and the flow tends to self-organise itself into large-scale structures; cf. [26].

A more quantitative comparison of the distribution of small- and large-scale structures can be obtained by considering the energy spectra; see Figure 9. A mutual characteristic of all simulations is that kinetic energy, which is initially concentrated in high wave numbers (small eddies), with time, is transferred to smaller wave numbers (large eddies). For small wave numbers, and in agreement with the literature [33], one can observe a behaviour $E(\kappa) \sim \kappa^3$ which is connected to spectral backscatter. After attaining a maximum energy, the spectra show a decaying behaviour from which, with decreasing viscosity (increasing Reynolds number), the classical $E(\kappa) \sim \kappa^{-3}$ slope can be determined. However, the slope is sometimes slightly steeper and thus shows more a $\kappa^{-4}$ behaviour at some later time instances. As mentioned in Section 4.2, a similar phenomenon of a slope between $\kappa^{-3}$ and $\kappa^{-4}$ has also been observed for the Kelvin–Helmholtz instability.

Returning to Figure 8, the last observation we want to make concerns the decaying nature of the flow. A smaller viscosity implies less molecular diffusion and, therefore, the colour bars show that the maximum and minimum values of the vorticity are considerably higher than for larger viscosities. In this context, Figure 7 shows the evolution of the kinetic energy and enstrophy over time for all three different simulations.

## 5 Summary and conclusions

In this work, we have considered inf-sup stable, exactly divergence-free $\boldsymbol{H}(\mathrm{div})$-conforming FEM for time-dependent incompressible flow problems. To the authors' knowledge, this represents the first contribution which examines $\boldsymbol{H}(\mathrm{div})$-FEM for *time-dependent viscous* incompressible flows. Our considerations can be split into a theoretical analysis part and an application-oriented numerical examples part.

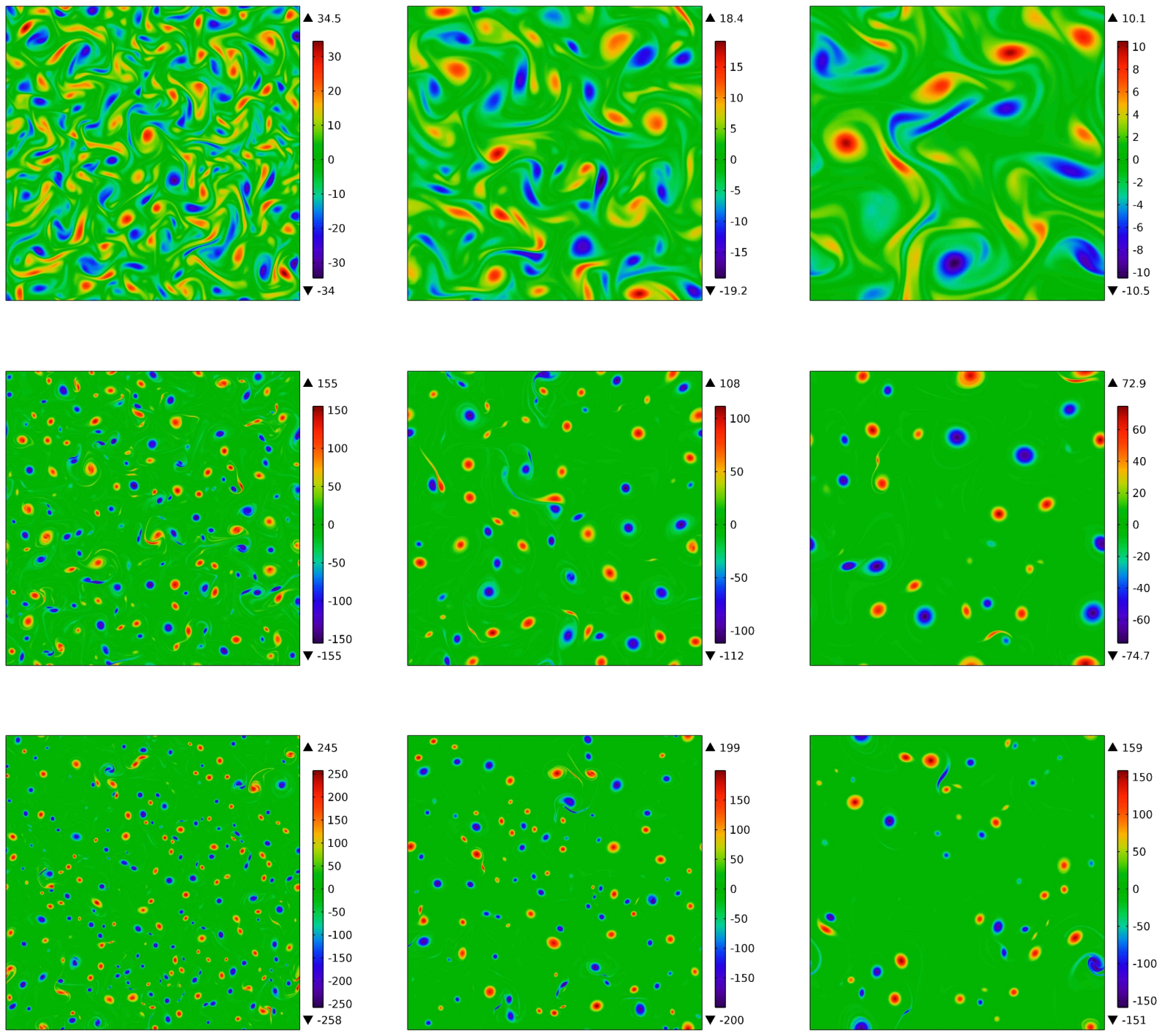

**Fig. 8** Vorticity $\nabla \times \boldsymbol{u}_h(t)$ for freely decaying 2D turbulence. From left to right: $t \in \{2,4,8\}$. From top to bottom: $\nu \in \{5 \times 10^{-5}, 10^{-5}, 4 \times 10^{-6}\}$. Results are obtained by div-free $\boldsymbol{H}(\mathrm{div})$-FEM RT3 with $\gamma = 1.5$ upwinding. Black triangles denote the maximum and minimum value attained over $\Omega$.

For the continuous-in-time Oseen problem, a numerical error analysis has been carried out with a special focus on obtaining pressure- and *Re*-semi-robustness. Because divergence-free methods allow for a separation of the velocity and pressure approximation, we only considered the velocity approximation in this work and, provided the exact solution is smooth enough, our derived estimates are of order $\mathcal{O}(h^k)$. A very important part of the analysis has been the usage of the discrete, stationary Stokes projection for the error splitting. For the Oseen problem it has been possible to show that the growth of the error w.r.t. time is only linear.

In the future, keeping pressure- and *Re*-semi-robustness in mind, we clearly intend to analyse the nonlinear Navier–Stokes problem, as well. In doing so, the discrete Helmholtz projection might be applied to fine-tune the error estimates. Furthermore, taking a closer look at the pressure approximation and discovering if divergence-free $\boldsymbol{H}(\mathrm{div})$-FEM also hold advantages for the pressure could be interesting.

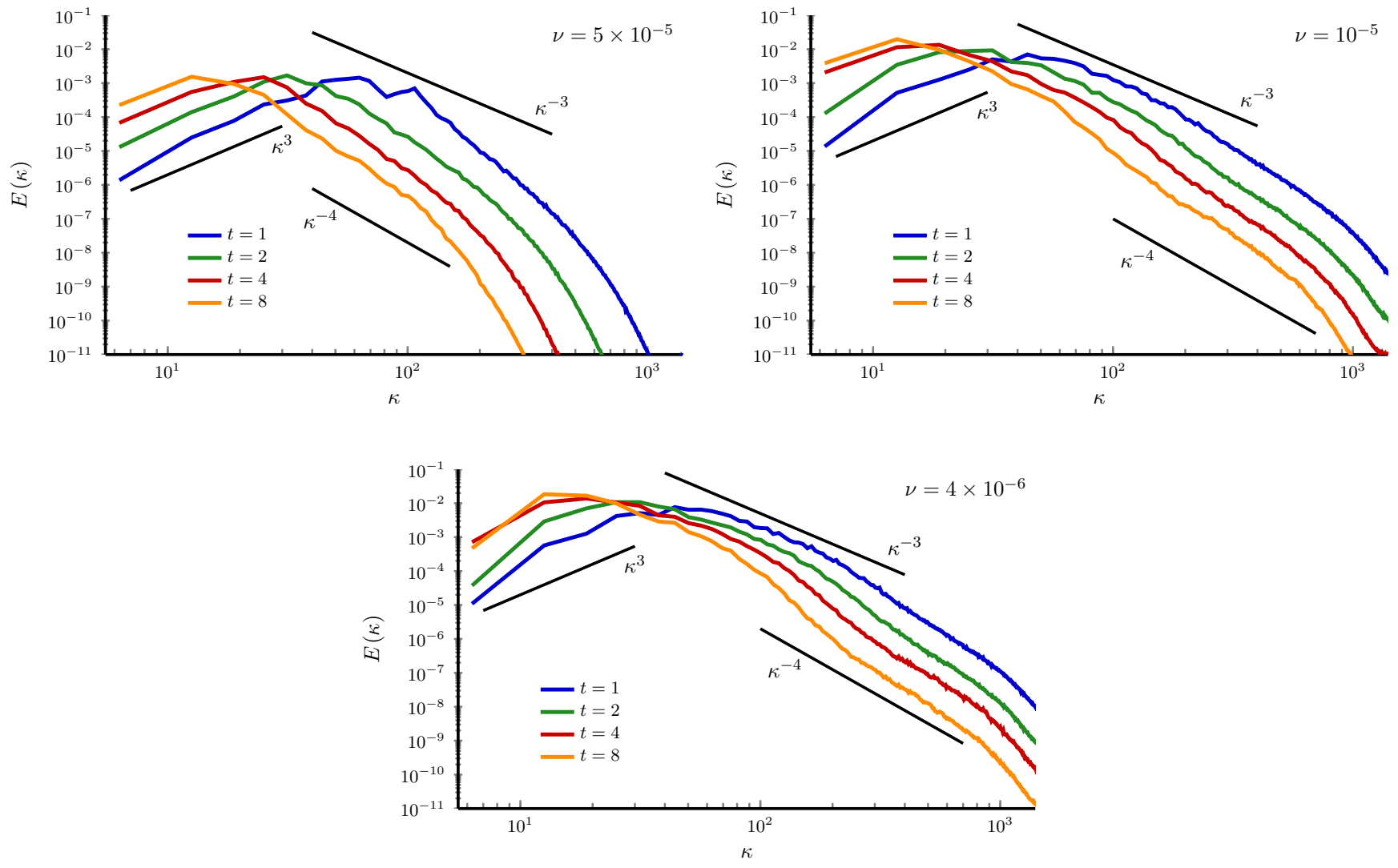

**Fig. 9** Energy spectrum $E(\kappa)$ and wave number $\kappa$ for $\nu \in \times\left\{5\times 10^{-5}, 10^{-5}, 4\times 10^{-6}\right\}$ at particular time instances. Spectra are computed with the MATLAB-based post-processing toolbox AnaFlame [50].

For the applications part, on the other hand, we have considered the nonlinear Navier–Stokes equations with a particular emphasis on two-dimensional high Reynolds number problems, which possess profound vortical structures. The first example, the planar lattice flow, has revealed that higher-order $\boldsymbol{H}(\mathrm{div})$-FEM do have advantages over more common finite element schemes. Especially the incorporation of a non-dissipative velocity jump upwind stabilisation for dominant convection seems to be very attractive and efficient. Contrary to the theory for the Oseen problem, this example has also shown that the error of the Navier–Stokes problem increases exponentially in time. Furthermore, $\boldsymbol{H}(\mathrm{div})$-FEM have been applied to the simulation of 2D Kelvin–Helmholtz instabilities, triggered by a plane mixing layer. We have shown that the problem is extremely sensitive in the sense that even though a very highly resolved method (3rd order FEM with nearly 5 million velocity DOFs) has been applied, such a thing as mesh convergence for the enstrophy is still not achieved and thus would be extremely expensive. However, the evolution of the kinetic energy is invariant with respect to mesh refinement. Lastly, the $\boldsymbol{H}(\mathrm{div})$-FEM have been applied to the simulation of freely decaying 2D turbulence. Our results are in agreement with theoretical considerations—both the behaviour of kinetic energy and enstrophy and the velocity spectra are consistent with previous research in this direction.

In the future, we intend to extend our numerical examples towards problems with no-slip conditions (in this work, we only considered periodic boundary conditions for the applications) and, therefore, problems involving boundary layers, separation and reattachment. In this context, also an extension to three-dimensional problems is planned where one has to deal with the aspect of efficient solvers. With regard to 3D problems, the question of suitable turbulence modelling also arises.

**Acknowledgements** We gratefully acknowledge the comments and suggestions about this paper from the anonymous reviewers.